\documentclass{article}
\usepackage{amsmath,amssymb,graphicx}
\newtheorem{thm}{Theorem}
\newtheorem{lem}[thm]{Lemma}
\newtheorem{lemma}[thm]{Lemma}

\newtheorem{prop}[thm]{Proposition}

\newcommand{\bi}{\begin{itemize}}
\newcommand{\ei}{\end{itemize}}
\newcommand{\be}{\begin{enumerate}}
\newcommand{\ee}{\end{enumerate}}
\newcommand{\bc}{\begin{center}}
\newcommand{\ec}{\end{center}}
\newcommand{\bt}{\begin{tabular}}
\newcommand{\et}{\end{tabular}}
\newcommand{\ba}{\begin{array}}
\newcommand{\ea}{\end{array}}
\newcommand{\ra}{\rightarrow}

\newcommand{\Z}{\mathbb Z}

 \newcommand{\specs}{$\mathrm{Spec}^{\mathrm{sum}}_k$}
 \newcommand{\specm}{$\mathrm{Spec}^{\mathrm{max}}_k(F) \ $}
 \newcommand{\spec}{$\mathrm{Spec}$}
  \newcommand{\sect}{\section}

  \newcommand{\bp}{\noindent\textit{Proof:}  }
  \newcommand{\ep}{\hfill$\Box$}

\newcommand{\BMS}{MR1929714}
\newcommand{\BurilloD}{MR1670622}
\newcommand{\josegrowth}{MR2102168}
\newcommand{\CannonFloydParry}{MR1426438}
\newcommand{\ProbList}{ProbList}
\newcommand{\CTClone}{MR2016187}
\newcommand{\FlajCF}{MR909335}
\newcommand{\Flajolet}{FlajBook}
\newcommand{\gubagrowth}{MR2104775}
\newcommand{\GubaSapirSubg}{MR1725439}

\newcommand{\Jit}{Jit}
\newcommand{\GenericCompl}{MR1981427}
\newcommand{\MiasBraid}{MR2423197}
\newcommand{\MiasArx}{MiasArx}
\newcommand{\gfun}{SaZi94}
\newcommand{\Stanley}{MR1676282}
\newcommand{\WZ}{WZ1985}
\newcommand{\Woodruff}{Woodruff}
\newcommand{\Arzh}{MR1794135}
\newcommand{\ArzhOlsh}{MR1445193}

\title{Random subgroups of Thompson's group
$F$\footnote{
The first, second and fourth authors
received support from a Bowdoin College Faculty Research
Award. The first author acknowledges support from a PSC-CUNY Research Award. The second author acknowledges the support of the Algebraic Cryptography Center at Stevens Institute of Technology, Hoboken New Jersey during the writing of this article.
 The third author thanks NSERC of Canada for financial support.
The fourth author acknowledges support from
NSF grant DMS-0604645. }}

\author{
Sean Cleary
 \\
\small Department of Mathematics,  \\[-0.8ex]
\small The City College of New York and the CUNY Graduate Center, New York, NY 10031\\[-0.8ex]
\small\texttt{cleary@sci.ccny.cuny.edu}\\[1.2ex]
Murray Elder\footnote{Corresponding author}\\
\small Mathematics, \\[-0.8ex]
\small The University of Queensland, Brisbane, Australia \\[-0.8ex]
\small\texttt{murrayelder@gmail.com}\\[1.2ex]
Andrew Rechnitzer\\
\small Department of Mathematics,  \\[-0.8ex]
\small University of British Columbia, Vancouver, Canada\\[-0.8ex]
\small\texttt{andrewr@math.ubc.ca}\\[1.2ex]
Jennifer Taback
\\
\small Department of Mathematics, \\[-0.8ex]
\small  Bowdoin College, Brunswick, ME 04011\\[-0.8ex]
\small\texttt{jtaback@bowdoin.edu}
}

\date{\small \today. 
 }

\begin{document}
\maketitle

\begin{abstract}
We consider random subgroups of Thompson's group $F$ with respect to
two natural stratifications of the set of all $k$ generator
subgroups. We find that the isomorphism classes of
subgroups which occur with positive density are not the same for the two stratifications.
We give the first known examples of {\em persistent} subgroups,
whose isomorphism classes occur with positive density within the set
of $k$-generator subgroups, for all sufficiently large $k$.
Additionally, Thompson's group provides the first example of a group
without a generic isomorphism class of subgroup.
Elements of $F$ are represented uniquely by reduced pairs of
finite rooted binary trees.
We compute the asymptotic growth rate
and a generating function for the number of reduced pairs of trees,
which we show is D-finite and not algebraic. We then use the asymptotic growth to prove our density results.
\end{abstract}

\sect{Introduction}\label{sec:intro}

We investigate the likelihood of randomly selecting a particular $k$-generator subgroup of Thompson's group $F$, up to isomorphism.  This is made precise through a notion of {\em asymptotic density}.  This in turn involves a choice of {\em stratification} of the set of $k$-tuples of elements, which we view as generating sets for the subgroups, into spheres of size $n$. Intuitively,
the density of an isomorphism class of subgroup with $k$ generators is the probability that a randomly selected $k$-generator subgroup is in the class.

A $k$-generator subgroup $H$ of a group $G$ is called {\em generic} among all $k$-generated subgroups if a randomly selected subgroup of $G$ with $k$ generators is isomorphic to $H$ with probability which is asymptotically one.  Previous results on asymptotic density of subgroups of particular groups, such as braid or free groups, have always found a generic type of subgroup for all $k$.  We find that Thompson's group $F$, with respect to each of two natural stratifications on the set of $k$-generator subgroups, does not possess a generic isomorphism class of subgroup for any $k$.  Additionally, for each stratification there are isomorphism classes of subgroups which are chosen at random with  small but positive probability among the set of all $k$-generated subgroups, for {\em any} sufficiently large $k$.  We call such subgroups {\em persistent}.  Lastly, we exhibit subgroups with positive density with respect to one stratification but not the other, illustrating that different natural notions of stratification can have dramatic effects on the forms of  randomly chosen subgroups.

The likelihood that a particular isomorphism class of subgroup of a given group is
selected at random is motivated by questions in group-based
cryptography. The analysis of the security of algorithms used in cryptography can depend upon the
expected isomorphism type of a random subgroup. Many group-based cryptosystems propose the braid group $B_n$ as a platform; recent work of Miasnikov, Shpilrain and Ushakov \cite{\MiasBraid} shows that experimentally, subgroups of $B_n$ generated by $k$ elements where $k$ is small relative to $n$, and moreover, those $k$ elements are of small size, are generically isomorphic to $B_n$.  Due to the restrictions on the size of the generators we cannot conclude that a subgroup of $B_n$ with $k$ generators is generically isomorphic to $B_n$.  Regardless, their results  explain why current  cryptosystems based on $B_n$ are  vulnerable to attack.

Our definition of the asymptotic density of a particular subgroup $H$ of a group $G$
follows Borovik, Miasnikov and Shpilrain in \cite{\BMS}. They present a detailed discussion of asymptotic and
statistical questions in group theory. We also refer the reader to
Kapovich, Miasnikov, Schupp and Shpilrain \cite{\GenericCompl} for
background on generic-case complexity and notions of density.

We let $G$ be an infinite group and $X$ a set of representatives of
elements that maps onto $G$. We can associate to each $x \in X$ an
integer {\em size}. For example,  a natural notion of size is word length- we can let $X$ be the set of all words in
a finite generating set for $G$, with size corresponding to word
length.  There are situations where other notions of size, besides word length, are considered.
 We let $X_k$ be the set of unordered $k$-tuples of
representatives $x \in X$. Then each member of $X_k$ corresponds to
a $k$-generated subgroup of $G$, taking the $k$ representatives as
the generators. We fix a notion of size on $X_k$.
We can define an integer size for each $k$-tuple in a variety of
ways. For example, the size of a $k$-tuple could be the sum of the
sizes of its components. Alternatively,  one could take the size of a
$k$-tuple to be the maximum size of any of its components. Once
notions of size are fixed, both for elements and tuples, the set of all tuples of size $n$ in $X_k$
is called the $n$-sphere, and denoted Sph$_k(n)$. Such a
decomposition of $X_k$ into spheres of increasing radii is known as
a {\em stratification} of $X_k$. We prefer our spheres of a fixed
size to be finite and thus we can regard these spheres of increasing
radii as an exhaustion of an infinite set $X_k$ by a collection of
finite sets.

To quantify the likelihood of randomly selecting a particular subset
of $X_k$, we take a limit of the counting measure on spheres of
increasing radii.  Let $|T|$ denotes the size of the set $T$. The {\em asymptotic density} of a subset $T$ in
$X_k$ is defined to be the limit
\[\lim_{n\ra\infty} \frac{|T\cap \mathrm{Sph}_k(n)|}{|\mathrm{Sph}_k(n)|}\]
if this limit exists.
We often omit the word asymptotic and refer to this limit simply as
the {\em density} of $T$.

To understand density not just of $k$-tuples, but of isomorphism classes of $k$-generator subgroups, we let $T_H$ be the set of $k$-tuples that generate a subgroup of $G$
isomorphic to some particular subgroup $H$. If the density of $T_H$
is positive we say that $H$ is {\em visible} in the space of
$k$-generated subgroups of $G$. We call the set of all visible
$k$-generated subgroups of $G$ the {\em $k$-subgroup spectrum},
denoted by \spec$_k(G)$. If the density of $T_H$ is one, we say that
$H$ is {\em generic} in \spec$_k(G)$; if this density is zero we say that $H$ is
{\em negligible} in \spec$_k(G)$.

We make a series of choices within this construction, each of
which can greatly influence the densities of different subsets;
those choices include: the representation of group elements, the
size function defined on $X$, and the stratification of the set of tuples
$X_k$. Additionally, we are asserting that the likelihood of
randomly selecting a $k$-generator subgroup isomorphic to the given
one is captured by the limit as defined. It is certainly possible to
construct contrived stratifications which various pathological
properties,
so  we
concentrate on stratifications which correspond to ``natural''
definitions of the sphere of size $n$ in $X_k$.  Despite this, we  show
that for Thompson's group $F$, a small change in the stratification
has a great impact on the set of visible subgroups.

Below, we show that  Thompson's group $F$ is
the first example of a group which has different asymptotic
properties with respect to two different, yet natural, methods of
stratification.  To define these stratifications, we represent
elements of $F$ using reduced pairs of finite rooted binary trees,
which we abbreviate to ``reduced tree pairs".  These representatives
are in one-to-one correspondence with group elements.
Each pair
consists of two finite, rooted binary trees with the same number of
leaves, or equivalently, with the same number of internal nodes or
carets, as defined below, satisfying a reduction condition specified in Section \ref{sec:comb}.
The size of a tree pair will be the number
of carets in either tree of the pair.

Using  reduced tree pairs to represent elements of $F$, we
define the sphere of radius $n$ in $X_k$ in two natural ways:
\begin{enumerate}
\item  take Sph$_k(n)$ to be the set of $k$-tuples in which the sum
of the sizes of the coordinates is $n$, or
\item take Sph$_k(n)$ to be the set of $k$-tuples where the maximum size of a coordinate is
$n$.
\end{enumerate}
We will refer to these as the ``sum stratification'' and ``max stratification'' respectively.
With respect to the sum stratification, every non-trivial
isomorphism class of $m$-generated subgroup  for $m\leq k$ is
visible. That is,  every possible subgroup isomorphism class has non-zero density.  With respect to the max
stratification, there are subgroup isomorphism classes with zero density.

Perhaps the most natural stratification to consider on $F$, or on any
finitely generated group, is obtained by taking the size of an
element of $F$ to be the word length with respect to a particular
set of generators.  For $F$ we can consider word length with respect
to the standard finite generating set $\{x_0,x_1\}$. This stratifies
the group itself into metric spheres.  Despite work of Jos\'e
Burillo \cite{\josegrowth} and Victor Guba \cite{\gubagrowth} in
this direction, the sizes of these spheres have not been calculated,
and thus it is not yet computationally feasible to consider the possible induced
stratifications of $X_n$ with respect to word length as a notion of size.

It is striking in our results below that the $k$-generator subgroups
of Thompson's group $F$ have no generic isomorphism type with
respect to either stratification, for any $k$. All other groups which have been
studied in this way exhibit a generic type of subgroup with respect
to  natural stratifications.  Arzhantseva and Olshanskii \cite{\ArzhOlsh} and Arzhantseva \cite{\Arzh}
considered generic properties of subgroups of free groups.  With respect to the notions
of stratification described here,
Jitsukawa \cite{\Jit} proved that $k$ elements of any finite rank free group
generically form a free basis for a free group of rank $k$.
Miasnikov and  Ushakov
\cite{\MiasArx} proved this is true also for the pure braid groups and right
angled Artin groups.

To obtain our results on random subgroups of Thompson's group $F$ we must be able to count the number $r_n$ of reduced pairs of trees with a given number of carets. Woodruff \cite{\Woodruff}, in his
thesis, conjectured that the number $r_n$ is proportional to $(8+4\sqrt 3)^n/n^3$. We prove
Woodruff's conjectured growth rate, and additionally show that the
generating function for the number of reduced tree pairs is not algebraic, but that it is D-finite, meaning that it satisfies a linear ordinary differential equation with polynomial coefficients.

This paper is organized as follows. In Section~\ref{sec:comb},  we
consider the number of pairs of reduced trees of size $n$, which we call $r_n$.
We prove that $r_n$ has a D-finite
generating function which is not algebraic. We prove that $r_n$
approaches $A\mu^n /n^3$ uniformly, where $A$ is a constant and $\mu
= 8 + 4 \sqrt{3} \approx  14.93$.

In Section~\ref{sec:F} we describe
particular subgroups of Thompson's group $F$ and elementary observations about $F$ that will be important in later
sections.

In Section~\ref{sec:sum} we study the sum stratification and compute the asymptotic density of
isomorphism classes of $k$-generator subgroups.  We prove that if $G$ is a non-trivial
$m$-generator subgroup of $F$, then its isomorphism class is visible in the
space of $k$-generator subgroups of $F$ for $k \geq m$.   This stands in stark contrast to
previously known examples, since no subgroup  is generic in this
stratification.

In Section~\ref{sec:max} we turn to the max stratification and compute the asymptotic density of
isomorphism classes of $k$-generator subgroups of $F$ and find very different behavior.  In this
case, not every isomorphism class of $m$-generator subgroup is visible  in the
space of $k$-generator subgroups of $F$ for $k \geq m$.  We
prove that $\Z$ is visible in the set of $k$-generated subgroups
only for $k=1$. Yet there are examples of isomorphism classes of
subgroups which are {\em persistent}; that is, visible in the set of
$k$ generator subgroups for all sufficiently large $k$.  For
example, we show that the isomorphism class of $F$ itself is visible
in the set of $k$-generated subgroups for all $k\geq 2$.

\textit{Acknowledgments:} The authors wish to thank  Collin Bleak,
Jos\'e Burillo, Jim Cannon, Steve Fisk, Bob Gilman, Alexei Miasnikov, Thomas Pietraho,
Claas R\"{o}ver, Mark Sapir, Melanie Stein   and Sasha Ushakov  for many helpful conversations
and feedback on this paper, and the anonymous referee for helpful suggestions.

\sect{Combinatorics of reduced tree pairs}\label{sec:comb}

\begin{figure}[ht!]
  \bc
 \texttt{a}\includegraphics[ scale=.35]{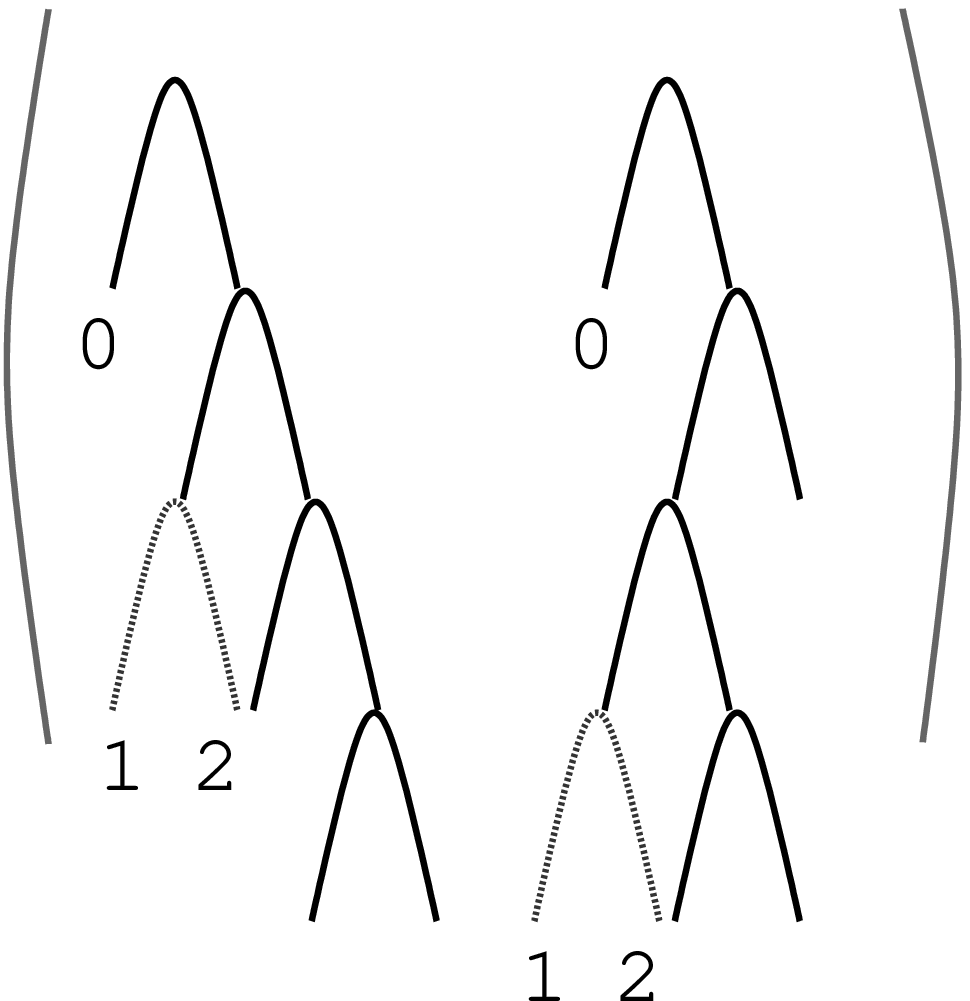}
\hspace{5mm} \texttt{b}\includegraphics[ scale=.35]{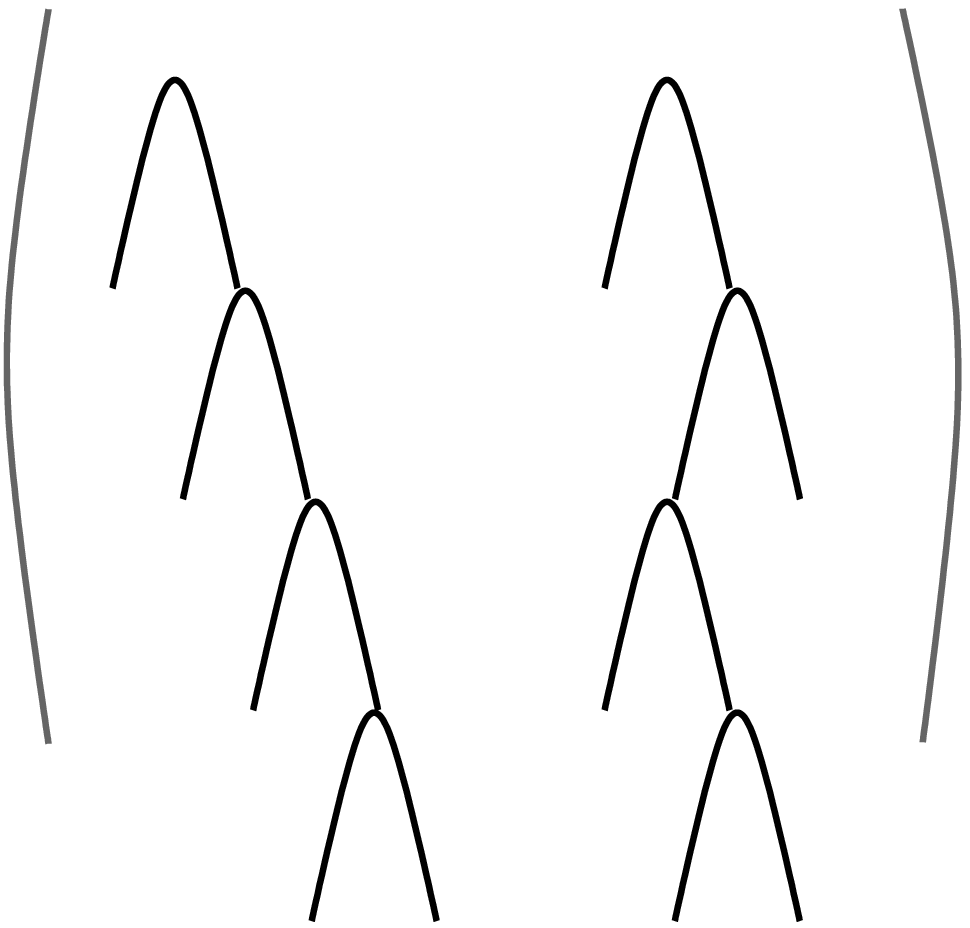}
\ec
  \caption{A five caret unreduced tree pair diagram, with its first three leaves numbered, and the corresponding four caret reduced tree pair diagram.}
  \label{fig:5caret}
\end{figure}

A {\em caret} is a pair of edges that join two vertices to a common
parent vertex, which we draw as $\wedge$.
 An $n$-caret {\em tree
pair diagram}, or {\em tree pair} for short, is an ordered pair
consisting of two rooted binary trees, each having $n$ carets.  A
5-caret tree pair is shown in Figure~\ref{fig:5caret}(a).
  A {\em leaf} is a vertex of degree one. A
tree with $n$ carets will have $n+1$ leaves.  In the trees we
consider, all vertices other than the leaves and the root have
degree three. The {\em left child} of a caret is the caret attached
to its left leaf; the {\em right child} is defined analogously.
An {\em exposed caret} is a caret both of whose children are leaves.
 A pair of trees with at least two carets in each tree is {\em unreduced} if, when
the leaves are numbered from left to right, each tree contains a
caret with leaves numbered $k$ and $k+1$ for some $k$. In an unreduced tree pair, the
caret with identical leaf numbers is removed from both trees, the leaves are
renumbered, and the trees are again inspected for possible reductions.
For example, the tree pair in  Figure~\ref{fig:5caret}(a)  is unreduced.
Removing the exposed caret with leaves labeled $1,2$ in each tree
yields the reduced tree pair in Figure~\ref{fig:5caret}(b).
A pair of trees which is not unreduced is called reduced.
Note that we do not reduce a pair of single carets: we insist that our tree pairs are always nonempty.
We denote the number of reduced tree pairs with $n$-carets by
$r_n$, so we have $r_0=0$ and $r_1=1$.

Ben Woodruff studied the enumeration of  $\{r_n\}$ in his thesis
\cite{\Woodruff} where he derived a formula for $r_n$ (which he
denoted $N_n$), proved an upper bound of  $(8+4\sqrt 3)^n\approx
14.93^n$ and conjectured an asymptotic growth rate of $(8+4\sqrt
3)^n/n^3$. We take a different approach to counting $r_n$ and derive
a recursive formula in terms of $c_n^2$, where $c_n= \frac{1}{n+1}
{2n \choose n}$ is the $n$-th Catalan number. Working in terms of
generating functions for $r_n$ and $c_n^2$, we obtain a finite-order
differential equation which leads to  a finite polynomial recurrence
for $r_n$. From this we are able to prove the growth rate
conjectured by Woodruff. The key to this section is to show that the
generating function for $r_n$ is closely related to that for $c_n^2$
and many of the properties of the
 generating function for $c_n^2$ are inherited by that of $r_n$.

\begin{figure}[ht!]
  \bc\bt{c|c|c|c|c}
 \includegraphics[ scale=.18]{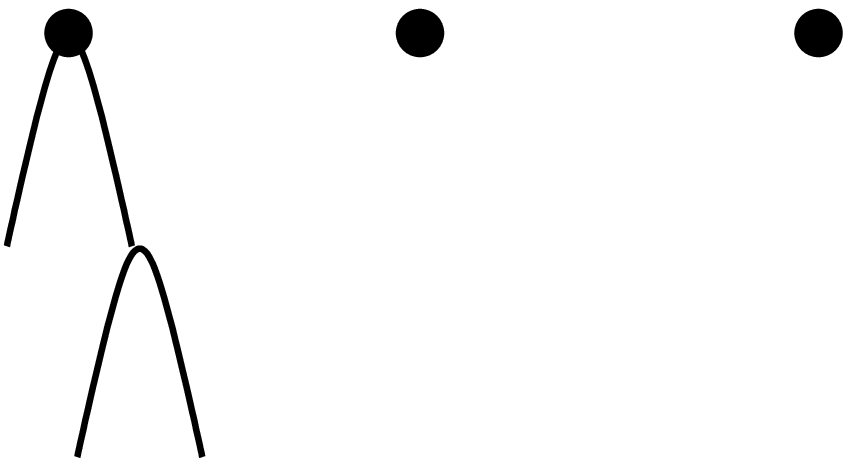}$\;\;$ & $\;\;$ \includegraphics[ scale=.18]{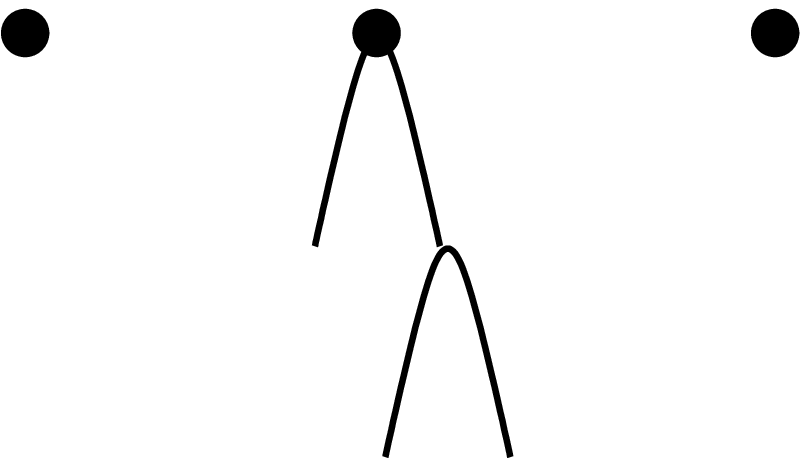}$\;\;$ & $\;\;$ \includegraphics[ scale=.18]{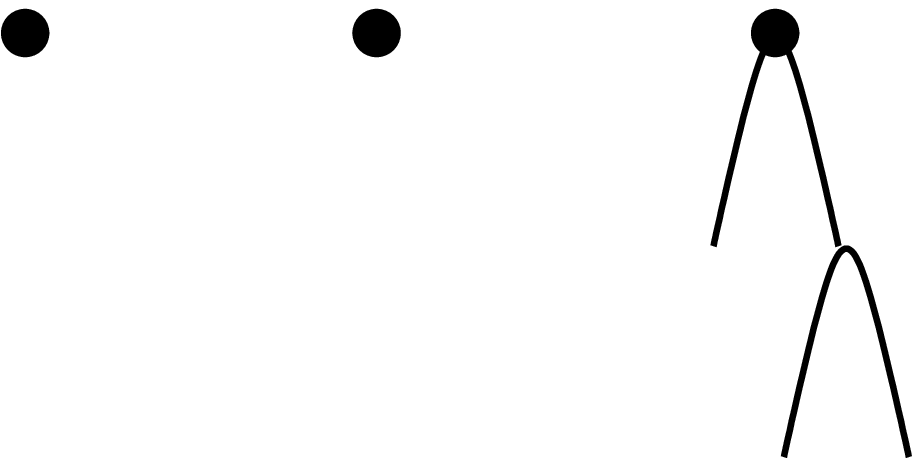}$\;\;$ & $\;\;$ \includegraphics[ scale=.18]{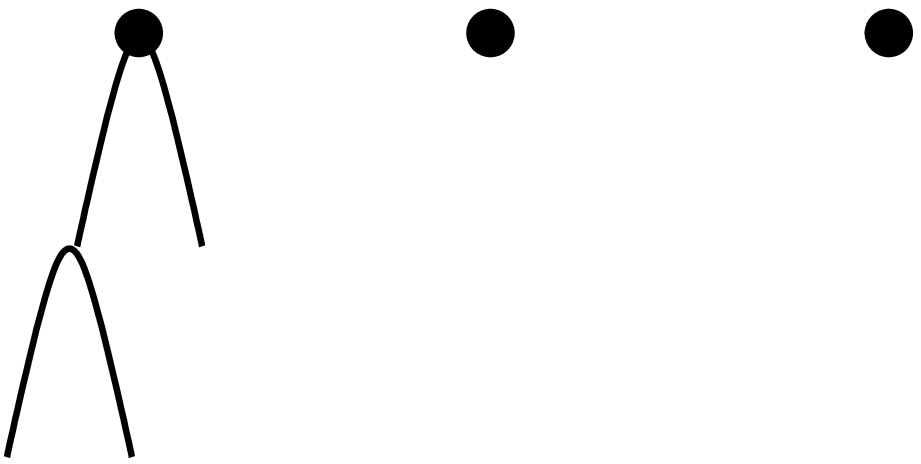}    $\;\;$ & $\;\;$ \includegraphics[ scale=.18]{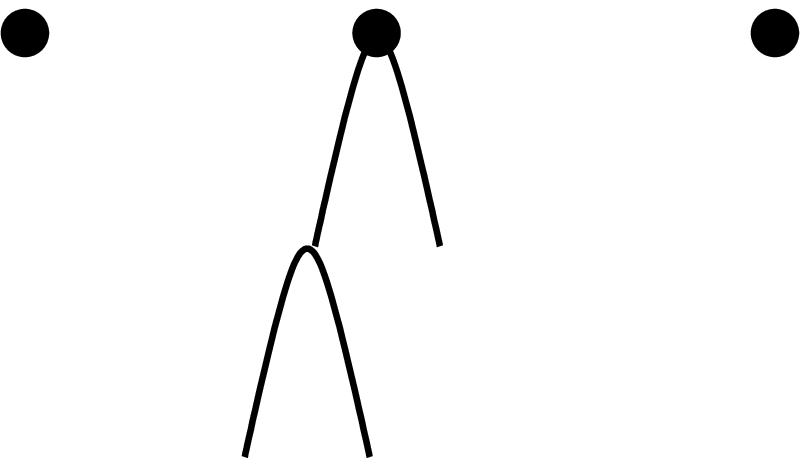}\et

  \vspace{5mm}

\bt{c|c|c|c}
 \includegraphics[ scale=.18]{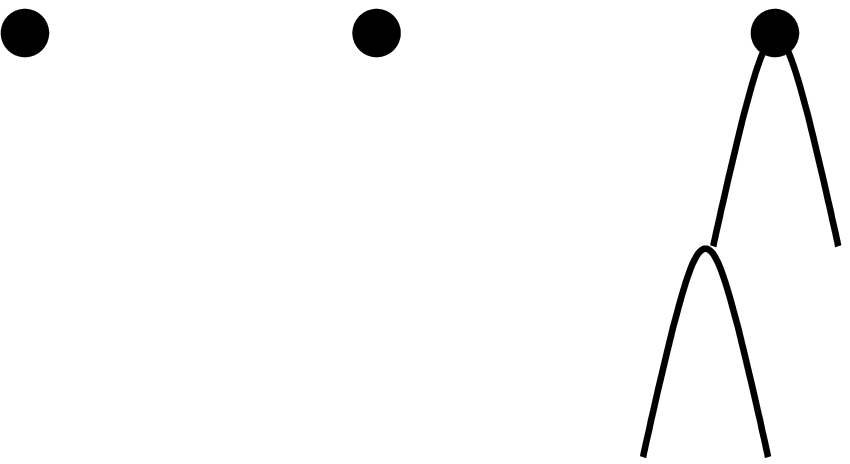}$\;\;$ & $\;\;$ \includegraphics[ scale=.18]{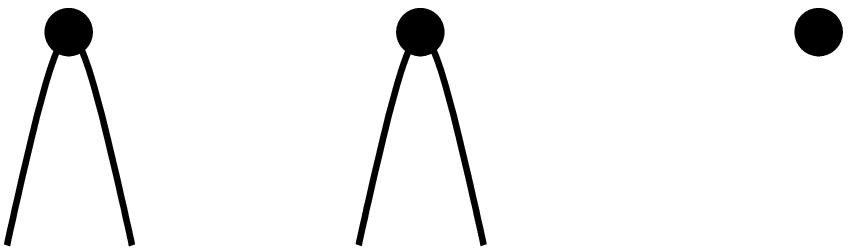}$\;\;$ & $\;\;$\includegraphics[ scale=.18]{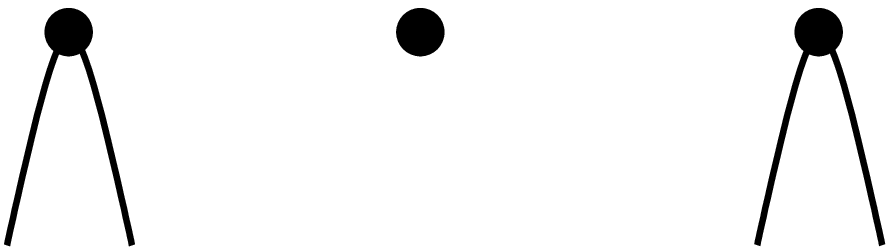}$\;\;$ & $\;\;$ \includegraphics[ scale=.18]{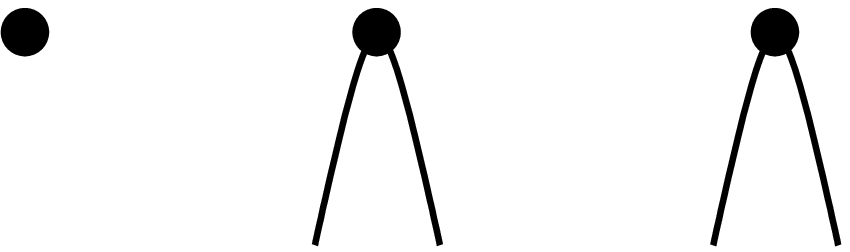}\et
\ec
  \caption{We exhibit that $f(3,2)=9$ by enumerating all forests consisting of three trees and a total of two carets.}
  \label{fig:forest32}
\end{figure}

We let $f(k,m)$ denote the number of ordered $k$-tuples of possibly empty rooted
 binary trees using a total of $m$ carets, which we call \emph{forests}. So for example $f(3,2)$, which is the number of forests of three trees  containing a total of two carets, is equal to nine, as
 shown by Figure~\ref{fig:forest32}. A straightforward argument shows that
$  f(k,n)=\frac{k}{2n+k}{2n+k\choose n}$.

\begin{figure}[ht!]
  \bc
 \includegraphics[ scale=.45]{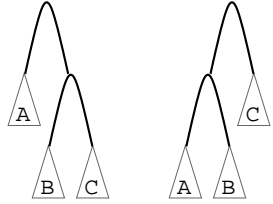}
\ec
  \caption{Decorating a reduced tree with a forest of three trees $A,B$ and $C$.}
  \label{fig:decorate}
\end{figure}

The $n$-th Catalan number $c_n$ counts the number of binary trees
consisting of $n$ carets, and thus $c_n^2$ is the number of ordered
pairs of rooted binary trees with $n$ carets in each tree. Some of
these pairs will be reduced, and some not. For those that are not
reduced, we can cancel corresponding pairs of carets to obtain an
underlying reduced tree pair. In a reduced tree pair consisting of $i$ carets,
each tree has $i+1$ leaves. We describe a process which is the
inverse of reduction, which we call ``decoration.''
To decorate a reduced
tree pair diagram $(S,T)$ with $i$ carets in each tree, we take a
forest of $i+1$ trees, some of which may be empty, and $n-i$ carets (for $n \geq i$), duplicate it, then append
the trees in the forests to the corresponding leaves of $S$ and $T$.
The first tree in the forest is appended to the first leaf, the
second tree in the forest to the second leaf and so on. We can do
this in $f(i+1,n-i)$ different ways. This decorating process yields
a new unreduced tree pair with $n$ carets, which will reduce to the
original reduced tree pair $(S,T)$ with $i$ carets. For example, the
reduced 2-caret tree pair drawn in bold in Figure~\ref{fig:decorate} can be
decorated in 9 different ways with a forest consisting of three trees $A$,$B$ and $C$ with a total of three carets between them, to yield unreduced pairs of 5 carets all of which would all reduce to the original tree pair diagram.  This leads to the following lemma.

\begin{lem}[Relating $r_n$ and $c_n^2$]\label{lem:recursive-r_n}
For $n\geq 1$
\begin{align*}
c_n^2=r_n+r_{n-1}f(n,1) +r_{n-2}f(n-1,2)+\ldots
+r_1f(2,n-1).\end{align*}
\end{lem}
\bp
Each $n$-caret tree pair is either reduced or must reduce to a unique
reduced tree pair of $i$ carets for some $i\in[1,n-1]$.
Hence the total number of  $n$-caret tree
pairs, $c_n^2$, is the number of pure reduced pairs of $n$-carets, $r_n$,
plus the number $r_i$ of reduced $i$-caret tree pairs multiplied by the
number of ways to decorate them with a forest of $n-i$ carets, $f(i+1,n-i)$, for
each possible value of $i$.
\ep

We can reformulate this recursion in terms of generating functions.
We define the generating functions for $r_n$, $c_n$ and $c_n^2$ respectively as:
\[\begin{array}{ll}
R(z)  = & r_1z+r_2z^2+r_3z^3+\ldots\\
C(z)  = &c_0+c_1z+c_2z^2+c_3z^3+\ldots\\
P(z) =  &c_1^2z+c_2^2z^2+c_3^2z^3\ldots
\end{array}\]
Note that $R(z)$ and $P(z)$ have no constant term while $C(z)$ does.
We prove in the following proposition that $R(z)$ can be obtained
from $P(z)$ via a simple substitution. Using knowledge of $P(z)$ we
can find a closed form expression for $R(z)$  and asymptotic growth
rate for $r_n$. Note that if $G(z)$ is the generating function for a
set of objects, then $G(z)^k$ is the generating function for ordered
$k$-tuples of those objects. In this way we can express the generating
function of $f(k,n)$ for fixed $k$ as $C(z)^k$.

\begin{prop}[Relating $R(z)$ and $P(z)$]
\label{prop:RzPzformula}
The generating functions for $r_n$ and $c_n^2$ are related by the following equation:
\begin{align*}
  R(z) &= (1-z) P( z(1-z) ),\\
  \intertext{which is equivalent to}
P(x) &= C(x) R(x C(x) ).
\end{align*}
\end{prop}
\bp
The generating function for the Catalan numbers is well known and may be written in closed form as
 \begin{align*}
C(x)=\frac{1-\sqrt{1-4x}}{2x};
\end{align*} it satisfies the algebraic equation
$C(x)(1-xC(x))=1$. See Stanley \cite{\Stanley} for example.
If we rewrite the equation
$R(z) = (1-z) P( z(1-z) )$
substituting the variable $z$ with $xC(x)$ then we obtain
\[R(xC(x))=(1-xC(x)) P\big( xC(x)(1-xC(x))\big)
=\frac1{C(x)}P(x)\] which rearranges to
\[P(x)=C(x)R(xC(x)) \]
This substitution is inverted by $x\mapsto z(1-z)$, and so proving
this  equation implies the proposition.  By
examining the coefficients of $x^n$ we will show that this statement
is equivalent to Lemma~\ref{lem:recursive-r_n}.

The right hand side can be written as \[C(x)R(xC(x)) =
C(x)\sum_{k=0}^\infty r_k \big(xC(x)\big)^k = \sum_{k=0}^\infty r_k
x^k \big(C(x)\big)^{k+1}\] We will use the notation $[x^i]G(x)$ to
denote the coefficient of $x^i$ in the expansion of a generating
function $G(x)$. Considering the above equation in terms of the
coefficient of $x^n$ we have
\begin{align}
c_n^2 = [x^n] P(x) &= [x^n] C(x) R(x C(x) )\nonumber\\
&= [x^n]   \sum_{k=0}^\infty r_k x^k \big(C(x)\big)^{k+1}
=  \sum_{k=0}^\infty [x^n] x^k r_k \big(C(x)\big)^{k+1}\nonumber\\
&= \sum_{k=0}^n [x^{n-k}]r_k  \big(C(x)\big)^{k+1}\nonumber
\end{align}

As noted above, $\big(C(x)\big)^{k+1}$ is the generating function
for the number of ordered $(k+1)$-tuples of rooted binary trees,
which are counted by $f(k+1,n)$. Thus the coefficient of $x^{n-k}$
in $\big(C(x)\big)^{k+1}$ is precisely $f(k+1,n-k)$, that is,
$[x^{n-k}]\big(C(x)\big)^{k+1}=f(k+1,n-k)$. So the above equation
becomes \[c_n^2=r_nf(n+1,0)+r_{n-1}f(n,1)+\ldots
r_1f(2,n-1)+r_0f(1,n)\] which is precisely
Lemma~\ref{lem:recursive-r_n} since $f(n+1,0)=1$ and  $r_0=0$.
\ep

A function is said to be {\em D-finite} if it satisfies a
homogeneous linear ordinary differential equation with polynomial
coefficients, for example, see \cite{\Stanley}.
The class of D-finite functions strictly contains the class of
algebraic (and rational) functions.  If one has a differential
equation for a generating function it is  possible to obtain
the asymptotic growth rate of its coefficients by studying the differential equation.  Following \cite{\Stanley}, a generating
function is D-finite if and only if its coefficients satisfy a
finite polynomial recurrence.

\begin{lem}[$R(z)$ is D-finite]
\label{lem:RDE}
The generating function $R(z)$ satisfies the following linear ordinary differential equation
\begin{multline*}
z^2(1-z)(16z^2-16z+1)(2z-1)^2\frac{\mathrm{d}^3R}{\mathrm{d}z^3} \\
-z(2z-1)(16z^2-16z+1)(8z^2-11z+5)\frac{\mathrm{d}^2R}{\mathrm{d}z^2} \\
-(128z^5-320z^4+365z^3-232z^2+76z-4)\frac{\mathrm{d}R}{\mathrm{d}z}\\
    +36z(z-1)R(z) = 0.
\end{multline*}
It follows that $R(z)$ is D-finite.
\end{lem}
\bp
Starting from a recurrence satisfied by the Catalan numbers we can
find a differential equation satisfied by $P(z)$ and then standard
tools allow us to transform this equation into one satisfied by
$R(z)$.

Since $c_n=\frac{1}{n+1}{2n \choose n}$, we have the following
recurrence for the Catalan numbers:
\begin{align*}
 (n+2) c_{n+1} &=  2(2n+1) c_n. \\
 \intertext{Squaring both sides yields}
 (n+2)^2 c_{n+1}^2 &= 4(2n+1)^2 c_n^2.
\end{align*}
Thus we have a finite polynomial recurrence for the coefficients of
$P(z)$, which means that we can find a linear differential equation for
$P(z)$. We do this using the Maple  package \texttt{GFUN} \cite{\gfun}
to obtain
\begin{align*}
(z^2-16z^3) \frac{\mathrm{d}^2P}{\mathrm{d}z^2} +
(3z-32z^2)\frac{\mathrm{d}P}{\mathrm{d}z} +(1-4z)P(z)
&=1.
\end{align*}
The original recurrence can be recovered by extracting the coefficient of $z^n$
in the above equation. We can then make this differential equation homogeneous
 \begin{align*}
(16z^3-z^2)\frac{\mathrm{d}^3P}{\mathrm{d}z^3}
+(80z^2-5z)\frac{\mathrm{d}^2P}{\mathrm{d}z^2}
+(68z-4)\frac{\mathrm{d}P}{\mathrm{d}z}
    +4P(z) &= 0.
\end{align*}

Making the substitution $z\mapsto z(1-z)$ using  the command
\texttt{algebraicsubs()} in \texttt{GFUN} we find a differential equation
satisfied by $P(z(1-z))$. This in turn leads to the homogeneous differential
equation for $R(z)$ given above.
\ep

Following the notation of Flajolet \cite{\Flajolet}, we say that two
functions are {\em asymptotically equivalent} and write $f(n)\sim
g(n)$ when \[\lim_{n\ra \infty} \frac{f(n)}{g(n)}=1.\]

\begin{prop}[Woodruff's conjecture]\label{prop:rate}
$r_n\sim A {\mu^n}/{n^3}$
where $\mu= 8+4\sqrt{3}$ and $A>0$ is a constant.
\end{prop}
\bp
We begin by establishing a rough bound on the exponential growth of
$r_n$ and refine this bound by analyzing a polynomial recurrence satisfied by
$r_n$ using techniques from Wimp and Zeilberger \cite{\WZ}.

Since reduced tree pairs are a subset of the set of all tree pairs, it follows
that $r_n \leq c_n^2$. We obtain a lower bound on $r_n$ by the following
construction. 
For each tree $T$  with $n$ carets, number the leaves from left to right starting with 0. Let $S_1$ denote the tree consisting of $n$ left carets, each the left child of its parent caret.  Let $S_2$ denote the tree with $n-1$ left carets, and a single interior caret attached to the right leaf of the leftmost caret.  This interior caret has leaves numbered 1 and 2. If $T$ does not have an exposed caret with leaves labeled 0 and 1, then the pair $(T,S_1)$  is reduced. If $T$ does have an exposed caret with leaves labeled 0 and 1, then form the reduced tree pair $(T,S_2)$.  Thus for each tree $T$ with $n$ carets, there is at least one distinct reduced tree pair diagram with $n$ carets, and we conclude that $c_n \leq r_n$.

It follows that $c_n^{1/n} \leq r_n^{1/n} \leq c_n^{2/n}$.
Since $c_n \sim B 4^n
n^{-3/2}$ for a constant $B$ (see Flajolet and Sedgewick \cite{\Flajolet} for example), it follows that $4 \leq \lim_{n
\to \infty} r_n^{1/n} \leq 16$.

The differential equation satisfied by $R(z)$ can be transformed into a linear
difference equation satisfied by $r_n$ using the Maple package \texttt{GFUN}
\cite{\gfun}:
\begin{align*}
 0 &= (n+5)(n+6)^2 r_{n+5} -(n+5)(n+4)(21n+101)r_{n+4} \nonumber\\
 &+2(4n+15)(n+4)(13n+33)r_{n+3} -4(n+3)(53n^2+208n+195)r_{n+2}\nonumber\\
 &+32(6n+5)(n+2)(n+1) r_{n+1} -64n^2(n+1)r_n.
\end{align*}

To compute the asymptotic behavior of the solutions of this recurrence we will
use the technique described in \cite{\WZ}.
This technique has also been automated by the command \texttt{Asy()} in the
\texttt{GuessHolo2} Maple package. This package is available from Doron
Zeilberger's website. We outline the method below.

Theorem~1 of \cite{\WZ} implies that the solutions of linear
difference equations
\begin{align*}
  \sum_{\ell=0}^\nu a(n) f_{n+\ell} &=0,
\end{align*}
where $a(n)$ are polynomials, have a standard asymptotic form. While this
general form is quite complicated (and we do not give it here),
 we note that in
the enumeration of combinatorial objects which grow exponentially rather
than super-exponentially one more frequently finds asymptotic expansions of
the form
\begin{align*}
 f_n &\sim \lambda^n n^\theta \sum_{j \geq 0} b_j n^{-j}.
\end{align*}
By substituting this asymptotic form into the recurrence one can determine
the constants $\lambda, \theta$ and $b_j$. For example,
substituting the above form into the recurrence satisfied by $r_n$, one obtains
(after simplifying):
\begin{align*}
 0 &= (\lambda-1)(\lambda^2-16\lambda+16)(\lambda-2)^2\nonumber\\
& +(\lambda-2)(5\lambda^4\theta+17\lambda^4-256\lambda^3-74\lambda^3\theta
+164\lambda^2\theta+558\lambda^2-352\lambda-96\lambda\theta+32)/n \nonumber\\
&+ O(1/n^2).
\end{align*}
In order to cancel the dominant term in this expansion we must have
\begin{align*}
 \lambda &= 1, 8 - 4\sqrt{3}, 2, 8 + 4\sqrt{3}.
\end{align*}
Each of these values for $\lambda$ implies different values of $\theta$ so as to cancel the
second-dominant term. In particular, if $\lambda=8+4\sqrt{3}$, then
$\theta=-3$. Since $4 \leq \lim_{n\ra \infty} r_n^{1/n} \leq 16$, it follows that the value of
$\lambda$ which corresponds to the dominant asymptotic growth of $r_n$ must be
$8+4\sqrt{3}$.

The application of this process using the full general asymptotic form
 has been automated by the \texttt{GuessHolo2} Maple package. In
particular, we have used the \texttt{Asy()} command to compute the asymptotic
growth of
$r_n$:
\begin{align*}
  \frac{n^3}{A \mu^n} r_n &\sim 1+ \frac{33/2 - 11\sqrt{3}}{n},
\end{align*}
for some constant $A$.
\ep

Though we do not need the exact value of the constant $A$ in
our applications below, we can estimate the constant $A$ as follows.
Using Stirling's approximation we know that $c_n^2 \sim \frac{1}{\pi n^3}
16^n$. This dictates the behavior of $P(z)$ around its dominant singularity, which forces the behavior of $R(z)$ around its dominant singularity.
Singularity analysis using methods of Flajolet and Sedgewick \cite{\Flajolet} then yields
\begin{align*}
 r_n &\sim \frac{6-3\sqrt{3}}{\pi n^3} \mu^n
 \sim \frac{12}{\mu \pi n^3} \mu^n.
\end{align*}
While this argument is not rigorous as it uses the estimate for $A$, the above form is in extremely close
numerical agreement with $r_n$ for $n \leq 1000$.

\begin{prop}[Not algebraic]
The generating function $R(z)$ is not algebraic.
\end{prop}
\bp
Theorem D of \cite{\FlajCF} states that if $l(z)$ is an algebraic
function which is analytic at the origin then its  Taylor
coefficients $l_n$ have an asymptotic equivalent of the form
\begin{align*}
l_n\sim A \, \beta^nn^s
\end{align*}
where $A \in \mathbb{R}$ and $s\not\in \{-1,-2,-3,\ldots\}$. Since $r_n$ is not of this form, in particular it has an $n^{-3}$ term, the generating function $R(z)$ cannot be
algebraic.
\ep

The generating function, or ``growth
series,'' for the actual word metric in Thompson's group $F$ with
respect to the $\{x_0,x_1\}$ generating set (see below), is not known to be algebraic
or even D-finite.
 Burillo \cite{\josegrowth} and Guba
\cite{\gubagrowth} have estimates for the growth but there are
significant gaps between the upper and lower bounds which prevent
effective asymptotic analysis at this time. Since finding
differential equations for generating
 functions can lead to information about the growth rate of the coefficients,
more precise understanding of the growth series for $F$ with respect the standard
 generating set (or any finite generating set) would be interesting and potentially quite useful.

In the following sections we regularly use following lemma which follows immediately from the
asymptotic formula for $r_n$.
\begin{lem}[Limits of quotients of $r_n$]
\label{lem:quotients}
For any $k\in \Z$
$$\lim_{n \rightarrow \infty} \frac{r_{n-k}}{r_n} = \mu^{-k}.$$
\end{lem}

\bp
From Proposition \ref{prop:rate} we have
\begin{align*}
r_{n-k}\sim A\mu^{n-k}(n-k)^3= A\mu^kn^3 \mu^{-k}\left(\frac{n-k}{n}\right)^3\sim r_n\mu^{-k}.
\end{align*}
\ep  


Finally, we give a formula for $r_n$. Woodruff (\cite{\Woodruff} Theorem 2.8) gave the following formula  for the number of reduced tree pairs on $n$ carets for $n\geq 2$
\begin{align*}
 \sum_{k=1}^{\lceil n/2\rceil}
2^{n-2k+1}{n-1\choose n-2k+1} c_{k-1}\sum_{ i=0}^k (-1)^i{k \choose i}c_{n-i}.
\end{align*}

One may readily verify (numerically) that Woodruff's formula and ours (below)
agree for $n \geq 2$. We have been able to show (using Maple) that both expressions
satisfy the same third-order linear recurrence, which together with the
equality of the first few terms is sufficient to prove that the expressions
are, in fact, equal. Unfortunately we have not been able to prove this more directly.

 \begin{lem}[Formula for $r_n$]\label{lem:formula-rn}
The number of reduced tree pairs with $n$ carets in each tree is
given by the formula
\begin{align*}
 r_n= \sum_{k= 1}^n  (-1)^{n-k}{k+1 \choose n-k}c_k^2
\end{align*}
 \end{lem}
 \bp
From Proposition~\ref{prop:RzPzformula} we have $R(z) = (1-z) P( z(1-z) )$ which expands to
\begin{align*}
\sum_{n\geq 1} r_nz^n &=(1-z) \sum_{k\geq 1} c_k^2z^k(1-z)^k\\
&= \sum_{k\geq 1} c_k^2z^k(1-z)^{k+1}\\
&=  \sum_{k\geq 1} c_k^2z^k  \left(\sum_{j=0}^{k+1}(-1)^j{k+1 \choose j}z^j\right)\\
&=  \sum_{k\geq 1} c_k^2\sum_{j=0}^{k+1}(-1)^j{k+1 \choose j}z^{k+j}
 \end{align*}
Now we look at the coefficient of $z^n$ on both sides. For the right
side, as $k$ runs from 1 up, we get exactly one term from the second
summation, when $j=n-k$. Thus we get
\begin{align*}
r_n &= \sum_{k\geq 1} c_k^2 (-1)^{n-k}{k+1 \choose n-k}
 \end{align*}
 which yields the result, since the binomial term becomes 0 for $k>n$.
\ep

\sect{Thompson's group $F$}\label{sec:F}

\begin{figure}[ht!]
  \bc \includegraphics[ scale=.27]{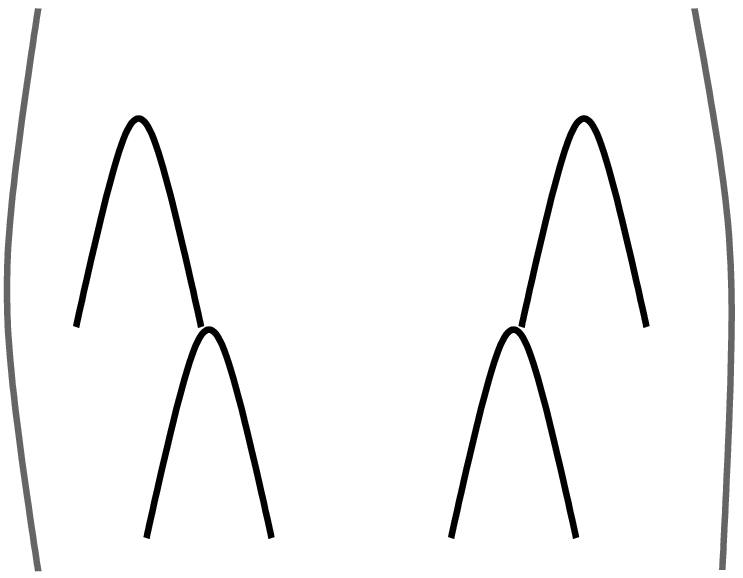} \hspace{2cm} \includegraphics[ scale=.27]{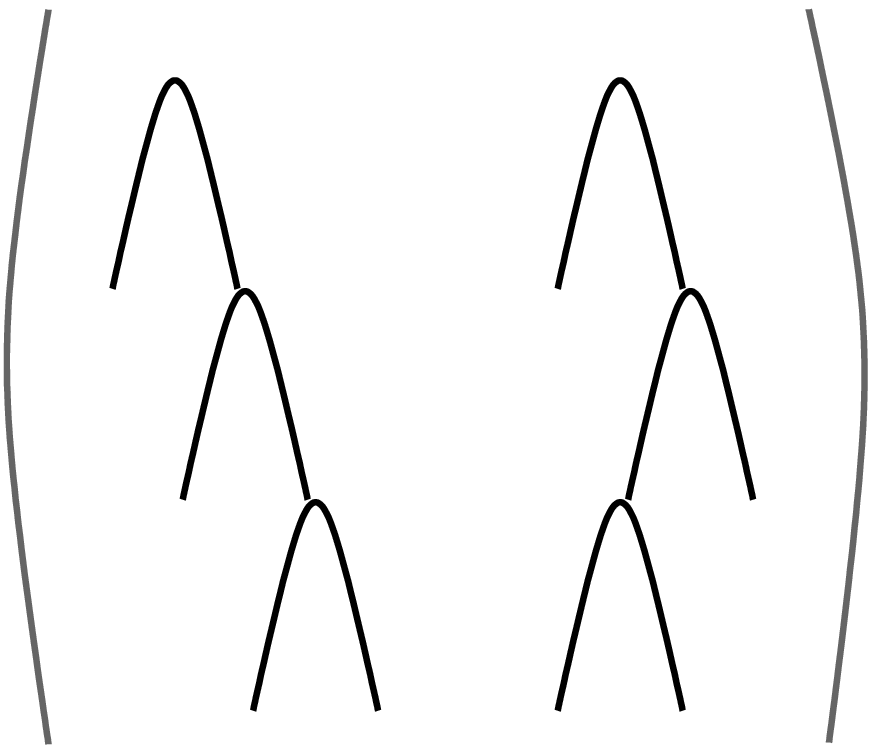}
      \ec
  \caption{Tree pair diagrams for the elements $x_0$ and $x_1$ respectively. }
  \label{fig:generators}
\end{figure}

Richard Thompson's group $F$ is a widely studied group which has
provided examples of and counterexamples to a variety of conjectures
in group theory.
We refer the reader to Cannon, Floyd and Parry \cite{\CannonFloydParry} for additional background information about this group.
Briefly, $F$ is defined using the standard infinite presentation
 \begin{align*}\langle x_0,x_1, \ldots | x_i^{-1}x_jx_i = x_{j+1}, \ i < j
 \rangle.\end{align*}
It is clear that $x_0$ and $x_1$ are sufficient to generate the entire group, and the standard finite presentation for this group is thus
\begin{align*}\langle
x_0,x_1|[x_0x_1^{-1},x_0^{-1}x_1x_0],[x_0x_1^{-1}, x_0^{-2}x_1x_0^2]
\rangle,
\end{align*}
where  $[a,b]$ denotes the commutator $aba^{-1}b^{-1}$. Group elements $w \in F$ can be uniquely represented by a reduced tree pairs as defined in the previous section. Equivalently, each element corresponds uniquely to a piecewise-linear map $\phi_w:[0,1] \rightarrow [0,1]$ whose slopes are all powers of two, the coordinates of the breakpoints are dyadic rationals and the slope changes at each breakpoint.  As described by Cannon, Floyd and Parry \cite{\CannonFloydParry}, each leaf of the reduced tree pair diagram defining $w \in F$ corresponds uniquely to an interval with dyadic endpoints in the domain or range of the map $\phi_w$.  The tree pair diagrams for $x_0$ and $x_1$ are given in Figure~\ref{fig:generators}.
$F$ has a diverse range of subgroups, but notably, it has no free subgroups of rank more than 1. 

\subsection{Recognizing support and commuting elements}

Two elements of $F$ can commute for many reasons, but one of the
 simplest is that they have disjoint supports.
The {\em support} of an element of $F$ regarded as a homeomorphism
of $[0,1]$ is the closure of the  set of points $x \in [0,1]$ such
that $f(x) \neq x$; that is, the set of points which are moved by
$f$.  Away from the support of $f$, the map $f$ will coincide with
the identity.  From the graph representing a group element as a homeomorphism,  it is easy to
recognize the complement of the full support of an element by inspecting where it
coincides with the identity;   $x_1$, for example, has support
$[1/2,1]$ as it coincides with the identity for the first half of
the interval.  It is not as easy
to recognize the complete support of an element directly from the reduced tree
pair diagrams representing it.  Nevertheless, it is possible to tell easily if the support
extends to the endpoints 0 and 1 of the interval, by inspecting
the locations of first and last leaves of the trees $S$ and $T$
representing an element.

 If the distances of the leftmost leaves (the leaves numbered 0) in $S$ and $T$ from
their respective roots are both $k$, then the homeomorphism represented by this
pair of trees coincides with
the identity at least on the interval $[0,\frac1{2^k}]$.
If there are, in addition to the leaves numbered 0, a sequence of
leaves numbered $1,\ldots,m$, each of which have the same distances from the root in both
trees, then the homeomorphism will coincide with the identity from
0 to the endpoint of the dyadic interval represented by leaf $m$.
 Similarly, near the right endpoint 1,
if the distances of the rightmost leaves (those numbered $n$) in $S$ and $T$ from
their respective roots are both $l$, then the homeomorphism represented coincides with
the identity at least on the interval $[1-\frac1{2^l}, 1]$.  Again, if there
are sequences of leaves numbered from $n-m$ up to $n$ which have the same levels in the trees $S$ and $T$, then the homeomorphism will coincide with the
identity on the corresponding dyadic interval, ending at the right endpoint of 1.
Elements that have homeomorphisms  that coincide with the identity for intervals of positive length at both the left and right endpoints are of particular interest as those elements lie
in the commutator subgroup of $F$, as described below.

\begin{figure}[ht!]
  \bc\bt{cccc}
    \includegraphics[ scale=.27]{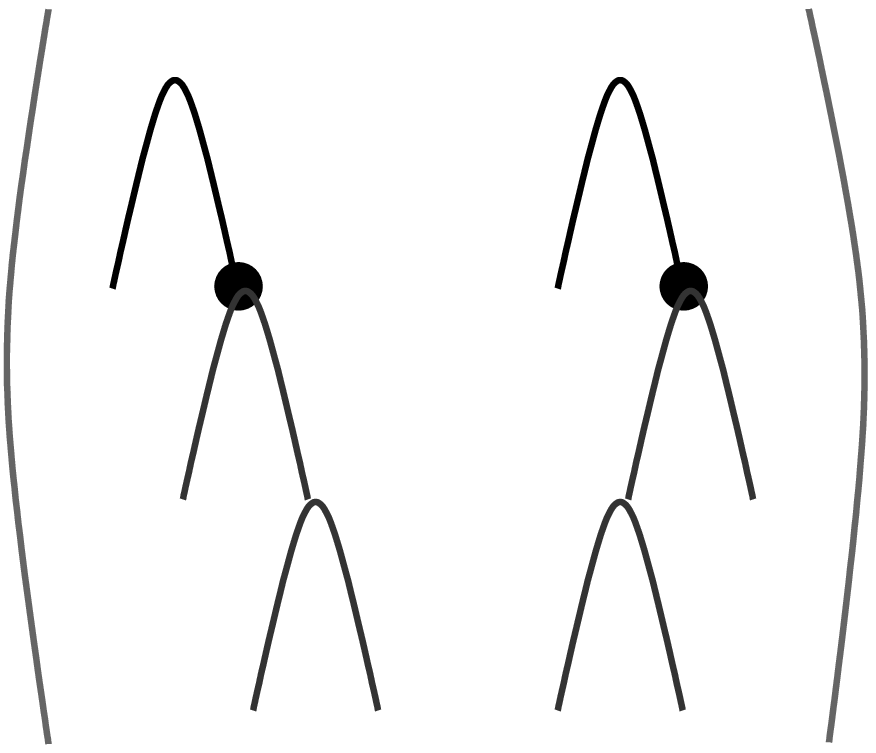}&
      \includegraphics[ scale=.27]{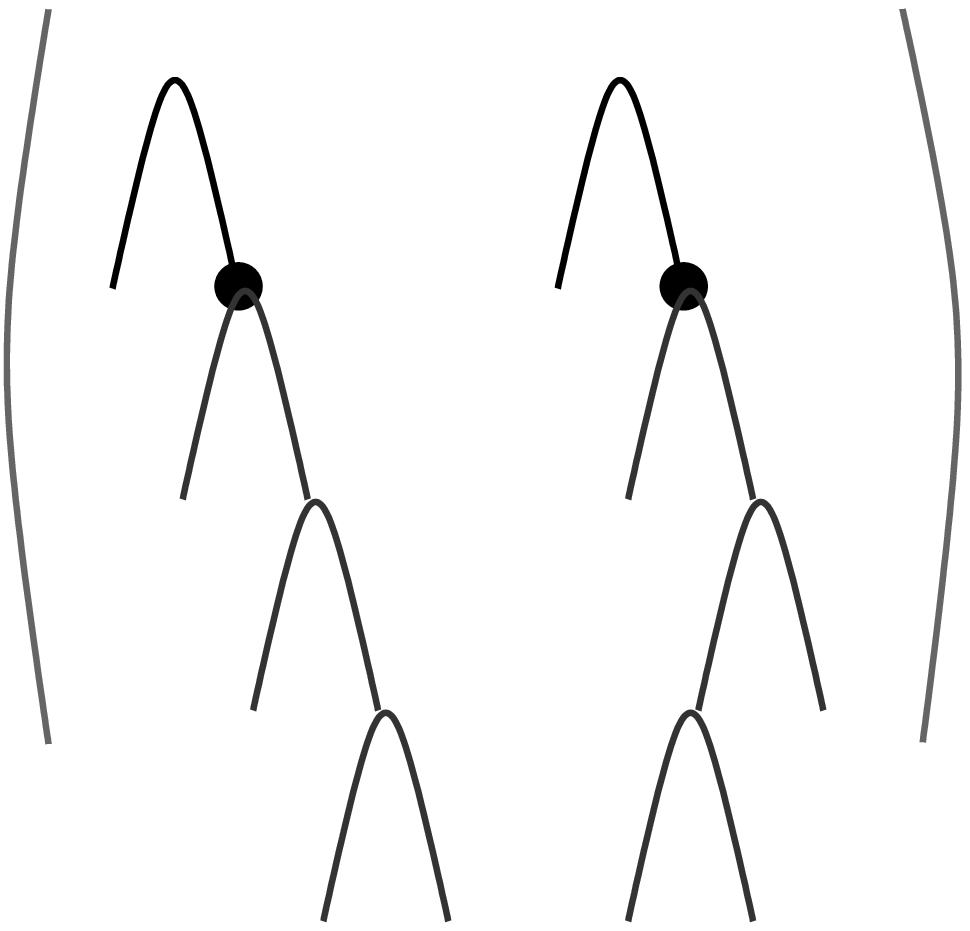}&
        \includegraphics[ scale=.27]{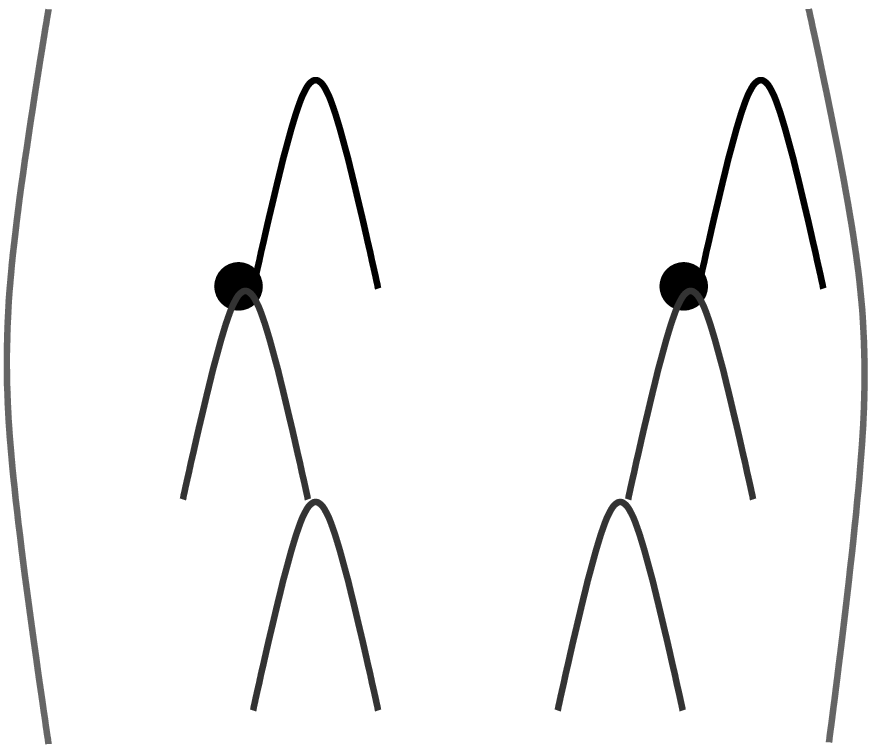}&
          \includegraphics[ scale=.27]{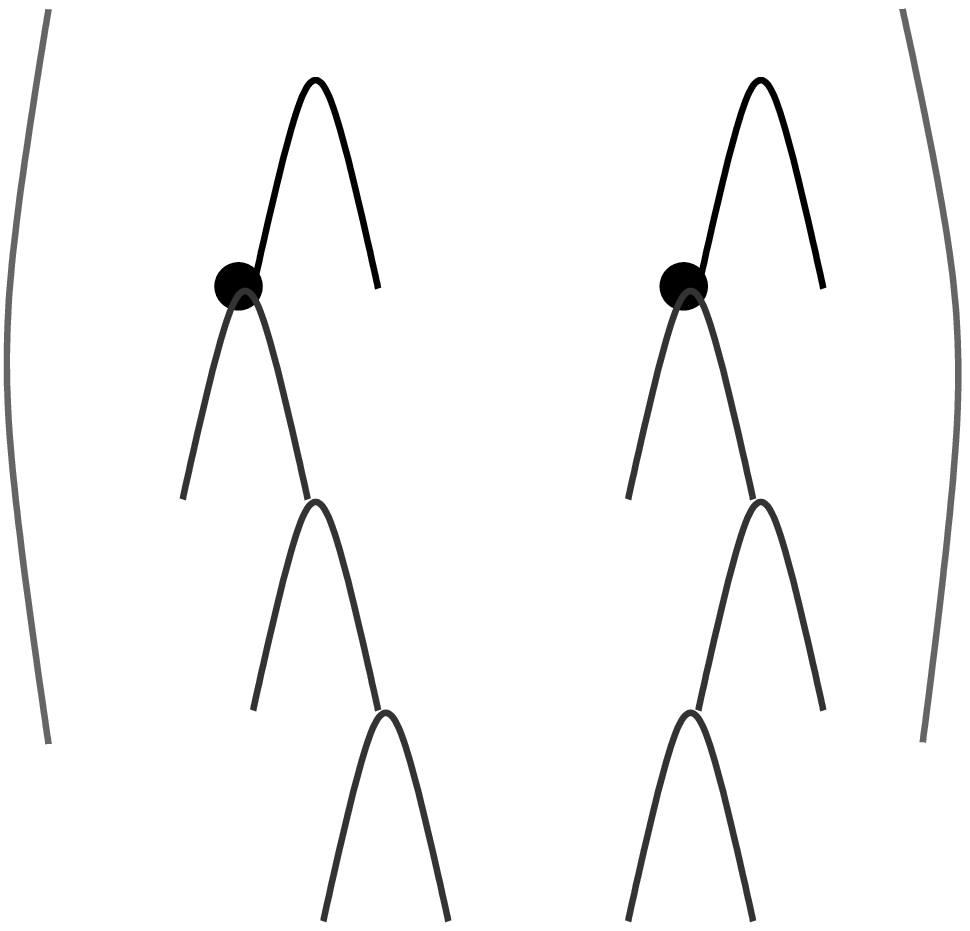} \\
       \et
      \ec
  \caption{Generators of the standard $F \times F$ subgroup of $F$.}
  \label{fig:FxF}
\end{figure}

 A simple method for generating pairs of commuting elements of $F$ is to construct them to have disjoint supports. An illustrative example is
simply the construction of a subgroup of $F$ isomorphic to $F \times F$,
where the four generators used are pictured in Figure \ref{fig:FxF}.
The first two generators have support lying in the interval
$[\frac12,1]$ and generate a copy of $F$ with support in that
interval.  Similarly, the second two
generators have support lying in $[0,\frac12]$ and generate a
commuting copy of $F$ in that interval. We refer to this example as the {\em
standard} $F \times F$ subgroup of $F$ and will make use of it in
later sections.

\subsection{More subgroups of $F$}
\label{sec:subgroups}

 One important subgroup of $F$ is the
restricted wreath product $\Z \wr \Z$. Guba and Sapir
\cite{\GubaSapirSubg} proved a dichotomy concerning subgroups of
$F$: any subgroup of $F$ is either free abelian or contains a
subgroup isomorphic to $\Z \wr \Z$. A representative example of a
subgroup of $F$ isomorphic to $\Z \wr \Z$ is easily seen to be
generated by the elements $x_0$ and $y = x_1x_2x_1^{-2}$.
The conjugates of $y$ by $x_0$ have disjoint support and thus commute.

Other wreath product subgroups of $F$ include $F \wr \Z$ and $H \wr
\Z$ for any $H < F$.  Generators for $H \wr \Z$ are obtained as
follows.  Let $\{h_1, \cdots ,h_k\}$ be a generating set for $H$
where $h_i = (T'_i,S'_i)$.  Let $T$ be the tree with two right
carets, and leaves numbered $1,2,3$.  Define generators
$k_i=(T_i,S_i)$ for $H \wr \Z$ by letting $T_i$ be the tree $T$ with
$T'_i$ attached to leaf $2$, and $S_i$ be the tree $T$ with $S'_i$
attached to leaf $2$.  Then $\{k_i\} \cup \{x_0\}$ forms a
generating set for $H \wr \Z$.

The group $F$ contains a multitude of subgroups isomorphic to $F$
itself; any two distinct generators from the infinite generating set
for $F$ will generate such a subgroup.  More generally, Cannon,
Floyd and Parry  \cite{\CannonFloydParry} describe a simple arithmetic condition to guarantee
that a set of analytic functions of the interval with the
appropriate properties generates a subgroup of $F$ which is
isomorphic to $F$.  A combinatorial description of their construction of proper subgroups of $F$ isomorphic to $F$ is as follows.

%

Given a finite string of zeros and ones, we construct a rooted binary tree by attaching to a root caret a left child if the first letter of the string is zero, and a right child otherwise.  Continue in this way, adding a child to the left leaf of the previous caret if the next letter in the string is a zero, to the right leaf of the previous caret otherwise.  For the final letter in the string, do not add a caret, but mark a distinguished leaf $v$ in the tree in the same manner, that is, mark the left leaf of the last caret added if the final letter is a zero, and the right leaf otherwise.  Let $T$ be a tree constructed in this way, and form two tree pair diagrams $h_0$ and $h_1$ based on $T$ as follows.  Denote $x_0 =
 (T_{x_0},S_{x_0})$ and $x_1 = (T_{x_1},S_{x_1})$. Draw four copies of
 the tree $T$, numbered $T_1$ through $T_4$.  To the marked
 vertex $v$ in $T_1$ attach the tree $T_{x_0}$ and to the marked vertex $v$ in
 $T_2$ attach $S_{x_0}$, forming the tree pair diagram
 representing $h_0$.  Do the same thing with $T_3, \ T_4, \
 T_{x_1}$ and $S_{x_1}$ respectively to form $h_1$.
Then $h_0$ and $h_1$ generate a subgroup of $F$ isomorphic to $F$,
which is called a {\em clone subgroup} in \cite{\CTClone} and
consists of elements whose support lies in the dyadic interval
determined by the vertex $v$. Subgroups of this form are easily seen
to be quasi-isometrically embedded.  This geometric idea is easily extended to construct subgroups of $F$
isomorphic to $F^n$.

\begin{figure}[ht!]
  \bc
    \includegraphics[ scale=.27]{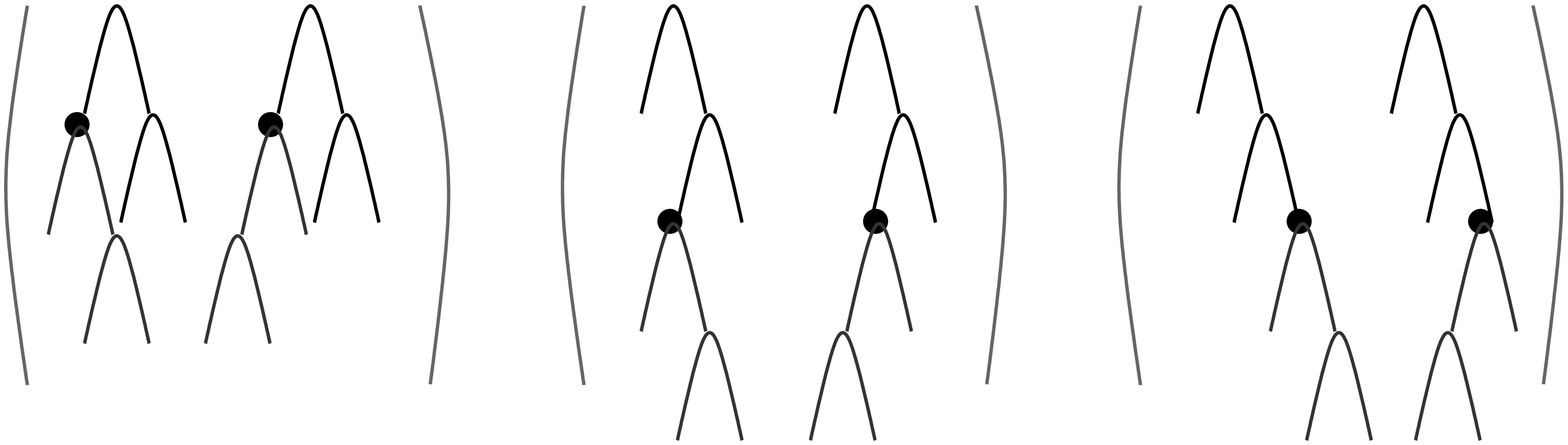}
      \ec
  \caption{Three tree pairs $h_1, h_2, h_3$ used to generate $\Z^3$.
  We have used the tree pair diagram for $x_0$ as each pair $(A_i,B_i)$. Note that the first
  pair can be reduced to a tree pair diagram containing only three carets by deleting the rightmost exposed caret.}
  \label{fig:Z3}
\end{figure}

Another  family of important subgroups of $F$ are the subgroups
isomorphic to $\Z^n$,
 which will play a role in the proofs in
Sections~\ref{sec:sum} and \ref{sec:max}.  We let $T$ be the tree
with $n-1$ right carets, and $n$ leaves, and $(A_i,B_i)$ for $i=1,2,
\ldots ,n$ reduced pairs of trees so that for each $i$, $A_i$ and $B_i$ have the same number of carets.  We construct generators $h_i=(C_i,D_i)$ of $\Z^n$ as
follows. We let $C_i$ be the tree $T$ with $A_i$ attached to leaf
$i$, and $D_i$ the tree $T$ with $B_i$ attached to leaf $i$, as
shown in Figure~\ref{fig:Z3}.
  We reduce the pair
$(C_i,D_i)$ if necessary.  It is easy to check by multiplying the
tree pair diagrams that $h_ih_j = h_jh_i$ for $i,j=1,2, \ldots ,n$
and thus these elements generate a subgroup of $F$ isomorphic to
$\Z^n$.  Burillo \cite{\BurilloD} exhibits a different family of
subgroups of $F$ isomorphic to $\Z^n$ using the generators
$\{x_0x_1^{-1}, x_2x_3^{-1},x_4x_5^{-1}, \ldots
x_{2n-2}x_{2n-1}^{-1} \}$ which he shows are quasi-isometrically
embedded.  In fact, Burillo proves that any infinite cyclic subgroup
of $F$ is undistorted; that is, that the cyclic subgroups are quasi-isometrically embedded.

\subsection{The commutator subgroup of $F$}\label{sec:commutator}
In the proofs in Sections~\ref{sec:sum} and \ref{sec:max} below, we
use both algebraic and geometric descriptions of the commutator
subgroup $[F,F]$.  This subgroup of $F$ has two equivalent
descriptions: \bi
\item The commutator subgroup of $F$ consists of all elements in $F$
which coincide with the identity map (and thus have slope 1) in neighborhoods
both of $0$ and of $1$. This is proven as
Theorem 4.1 of \cite{\CannonFloydParry}.
\item The commutator subgroup of $F$ is exactly the kernel of the
map $\varphi: F \rightarrow \Z \bigoplus \Z$ given by taking the
exponent sum of all instances of $x_0$ in a word representing $w \in
F$ as the first coordinate, and the exponent sum of all instances of
$x_1$ as the second coordinate.
\ei

The exponent-sum homomorphism $\varphi$ is closely tied to another
natural homomorphism $\phi$ from $F$ to $ \Z \bigoplus \Z$. The
``slope at the endpoints" homomorphism $\phi$ for an element $f \in
F$ takes the first coordinate of the image to be the logarithm base
2 of the slope of $f$ at the left endpoint 0 of the unit interval
and the second coordinate to be the logarithm base 2 of the slope at
the right endpoint 1. The images of the generators under the
slope-at-the-endpoints homomorphism $\phi$ are $\phi(x_0)= (1,-1)$
and $\phi(x_1)=(0,-1)$ and $\phi$ and
$\varphi$ have the same kernel.

\begin{figure}[ht!]
  \bc
\includegraphics[ scale=.45]{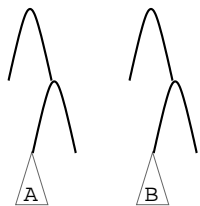}
\ec
  \caption{Constructing a tree pair representing a group element which
  lies in the commutator subgroup $[F,F]$.}
  \label{fig:commutator}
\end{figure}

It is not hard to see that the first description above has the
following geometric interpretation in terms of tree pair diagrams.
An  element of the commutator subgroup will have slope 1 at the left
and right endpoints and coincide with the identity on intervals of
the form $[0,b_0]$ and $[b_1,1]$ where $b_0$ and $b_1$ are,
respectively, the first and last points of non-differentiability in
$[0,1]$.  These points must lie on the line $y=x$, and the element
is represented by tree pair diagrams in which the first leaves
(numbered 0) in each tree lie at the same {\em level} or distance
from the root, and the same must be true of the last leaf in each of
the trees.  Thus, elements of the commutator subgroup are exactly
those which have a reduced tree pair diagram $(S,T)$ where the
leaves numbered zero are at the same level in both $S$ and $T$ and
the last leaves are also at the same level in both $S$ and $T$. For
example, if  $(A,B)$ is any reduced $n$-caret tree pair, then the
$(n+2)$-caret tree pair in Figure~\ref{fig:commutator} is also
reduced and represents an element in $[F,F]$.

We refer the reader to \cite{\CannonFloydParry} for a proof that the
commutator of $F$ is a simple group, and that $F/[F,F] \cong \Z
\bigoplus \Z$.

In our arguments below we will be interested in isomorphism classes
of subgroups of $F$.  It will sometimes be necessary to assume that
a particular finitely generated subgroup of $F$ is not contained in
the commutator subgroup $[F,F]$.  We now show that within the
isomorphism class of any subgroup $H$ of $F$, it is always possible
to pick a representative not contained in $[F,F]$.  The proof of this lemma follows the
proof of Lemma 4.4 of \cite{\CannonFloydParry}.

\begin{lem}[Finding subgroups outside the commutator]
\label{lem:outofcommutator} Let $H$ be a finitely generated subgroup
of $F$. Then there is a subgroup $H'$ of $F$ which is isomorphic to
$H$ and not contained in the commutator subgroup.
\end{lem}

\bp If $H$ is not contained in the commutator subgroup
$[F,F]$, then take $H' = H$.  Otherwise, let $H$ be generated by
$h_1,h_2, \ldots ,h_k$ where each $h_i \in [F,F]$.  Then each $h_i$
has an associated ordered pair $(a_i,b_i)$ where $a_i$ is
$x$-coordinate of the first point of non-differentiability of $h_i$
as a homeomorphism of $[0,1]$ (necessarily at $a_i$ the slope will
change from 1 to something which is not 1.) Similarly, we let $b_i$
be the $x$-coordinate of the final point of non-differentiability of
$h_i$. We let $a = \min \{a_i\}$ and $b = \max \{b_i \}$.  By the
choice of $a$ and $b$, all $h \in H$ have support in $[a,b]$.

Following the proof of Lemma 4.4 of \cite{\CannonFloydParry}, we let
$\phi:[a,b] \rightarrow [0,b-a]$ be defined by $\phi(x) = x-a$.  We
use $\phi$ to define a map on $h \in H$ by $h \mapsto \phi h
\phi^{-1}$, assuming that $\phi h \phi^{-1}$ acts as the identity
for $x \in (b-a,1]$.  It is clear from the definition of $\phi$ that
the breakpoints of $\phi h \phi^{-1}$ are again dyadic rationals,
and the slopes are again powers of two.  Since $\phi$ is an
isomorphism, we know that $H \cong \langle \phi h_i
\phi^{-1}\rangle$.  But this subgroup cannot be in the commutator,
since at least one element, the one which had its minimal breakpoint
at $x = a$, now has slope not equal to $1$ at $x=0$,
and thus is not in the commutator subgroup.
\ep

In the proofs in Sections~\ref{sec:sum} and \ref{sec:max} below, we
often want to make a more specific choice of representative subgroup
from an isomorphism class of a particular subgroup of $F$, as
follows.

Let $E_i(w)$ for $i=0,1$ denote the exponent sum of all instances of
$x_i$ in a word $w$ in $x_0$ and $x_1$.

\begin{lem}
\label{lem:exponentsum}
Let $H = \langle h'_1,h'_2, \ldots ,h'_k
\rangle$ be a finitely generated subgroup of $F$.  Then there is a
subgroup $H' = \langle h_1,h_2, \ldots ,h_k \rangle$ isomorphic to
$H$ so that $E_0(h_1) \neq 0$ and $E_0(h_j) = 0$ for $j=2,3, \ldots
,k$.
\end{lem}

\bp
By Lemma~\ref{lem:outofcommutator}, we assume without loss of
generality that $H$ is not contained in the commutator subgroup
$[F,F]$.  By replacing some generators with their inverses, we may
assume that $E_0(h'_i) \geq 0$ for all $i$, and that $E_0(h_1')$ is
minimal among those $E_0(h'_i)$ which are positive. For these
$h'_i$ with $i > 1$, we replace $h'_i$ by $h'_ih_1^{-d_i}$ where $d_i$ is chosen
so that $E_0(h'_ih_1'^{-d_i})$ is as small as possible while
non-negative. Repeating this process yields a generating set for a
subgroup isomorphic to $H$ with one element having exponent sum on
all instances of $x_0$ equal to zero. We can repeat this process
with the remaining generators, possibly reindexing at each step, until a generating set with the
desired property is obtained.
\ep

\sect{Subgroup spectrum with respect to the sum
stratification}\label{sec:sum}

We now introduce the first of  two stratifications of the set of $k$ generator
subgroups of Thompson's group $F$.  We view group elements as
non-empty reduced tree pairs  and denote by $X_k$  the
set of unordered $k$-tuples of non-empty reduced tree pairs $t_i = (T^i_1,T^i_2)$
for $i=1, \ldots ,k$. We denote the number of carets in $T^i_1$ by
$|t_i|$.  We define the sphere of radius $n$ in $X_k$ as the set of $k$-tuples
having a total of $n$ carets in the $k$ tree pair diagrams in the tuple:

\[\mathrm{Sph}_k^{\mathrm{sum}}(n)=\left\{(t_1,\ldots, t_k) \; | \; \sum_{i=1}^k |t_i|=n\right\}\]
which induces a stratification on $X_k$ that we will call the {\em sum stratification}.
 Note that since each tree in a tree-pair has the same number of carets, we only count (without loss of generality) the carets in the left tree.
For example, the triple
of tree pairs in Figure~\ref{fig:Z3}, once $h_1$ is reduced, lies in
$\mathrm{Sph}_3^{\mathrm{sum}}(11)$.

Recall from Section~\ref{sec:intro} that
the density of  a set $T$ of $k$-tuples of reduced tree pairs  is given by
\[\lim_{n\ra\infty} \frac{|T\cap \mathrm{Sph}_k^{\mathrm{sum}}(n)|}
{|\mathrm{Sph}_k^{\mathrm{sum}}(n)|}\] with respect to this
stratification. Let $H$ is a subgroup of $F$,  and $T_H$ the set of
$k$-tuples whose coordinates generate a subgroup of $F$ that is
isomorphic to $H$. Recall that $H$ is {\em visible} if $T_H$ has
positive density, and the {\em $k$-spectrum} \specs$(F)$ is the set
of visible subgroups with respect to the sum stratification of
$X_k$. In this section we explicitly compute these subgroup spectra.
We find that any isomorphism class of nontrivial subgroup $H$ of $F$ which can be
generated by $m$ generators is an element of in \specs$(F)$ for all $k \geq
m$ (Theorem~\ref{thm:specs}).  We conclude
 that this stratification does not distinguish any particular subgroups through
the subgroup spectrum, in contrast to the results we will describe
in Section~\ref{sec:max} when the max stratification is used.

We begin by determining upper and lower bounds on the size of the
sphere of radius $n$ in this stratification. Since
our $k$-tuples are unordered, we may assume that they are arranged from
largest to smallest.

\begin{lem}[Size of $\mathrm{Sph}^{\mathrm{sum}}_k(n)$]\label{lem:SizeSum}
For $k\geq 1$ and $n\geq k$, the size of the sphere of radius $n$
with respect to the sum stratification satisfies the following
bounds:
\[r_{n-k+1}\leq |\mathrm{Sph}^{\mathrm{sum}}_k(n)| \leq r_{n+k-1}.\]
\end{lem}
\bp
For the lower bound, $\mathrm{Sph}^{\mathrm{sum}}_k(n)$ contains all
$k$-tuples where the first pair has $n-k+1$ carets and the
remaining $(k-1)$ pairs are consist of two single carets. There are
$r_{n-k+1}$ ways to choose this first pair, which yields the lower
bound.

\begin{figure}[ht!]
  \bc
    \includegraphics[ scale=.27]{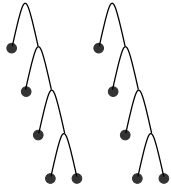}
      \ec
  \caption{A tree-pair consisting of $k-1$ right caret pairs (with $k=5$).}
  \label{fig:upperBd}
\end{figure}

For the upper bound, we consider the set of all $r_{n+k-1}$ reduced tree
pairs  with $n+k-1$ carets in each tree. A (small) subset of these
correspond to the $k$-tuples of $\mathrm{Sph}^{\mathrm{sum}}_k(n)$
as follows. Take the subset of these tree pairs where each tree
contains at least $k-1$ right carets, as in
Figure~\ref{fig:upperBd}, where leaf $i$ for $0 \leq i \leq n-1$ has
a possibly empty left subtree labeled $A_i$ in $T_-$ and $B_i$ in
$T_+$.  Let $A_n$ and $B_n$ respectively denote the right subtrees
attached to leaf $n$ in $T_-$ and $T_+$.  The sum of the number of
carets in the $A_i$ must equal $n$.

When the number of carets in $A_i$ equals the number of carets in
$B_i$ for all $i$, this pair of trees can be associated to an
(ordered) $k$-tuple of tree pairs with a total of $n$ carets.
Amongst these we can find every unordered $k$-tuple in
$\mathrm{Sph}^{\mathrm{sum}}_k(n)$. So this is a gross overcount
which suffices to prove the lemma.
\ep

\begin{thm}[All subgroup types are visible with respect to sum]
\label{thm:specs} Let \\ $H = \langle h_1,h_2, \ldots ,h_m \rangle$ be
a nontrivial subgroup of $F$.  Then $H \in$ \specs$(F)$ for all $k
\geq m$.
\end{thm}

We use the notation from Section~\ref{sec:commutator} to represent
the exponent sum of different generators in a word in $x_0$ and
$x_1$.  Let $E_i(w)$ for $i=0,1$ denote the exponent sum of $x_i$ in
a group element given by a word $w$.

\bp
Applying Lemmas~\ref{lem:outofcommutator} and
\ref{lem:exponentsum}, we may assume that $H$ is a representative
of its isomorphism class which is not contained in the commutator
subgroup $[F,F]$ and such that $E_0(h_1) \neq 0$ but $E_0(h_i) = 0$
for $i > 1$.

We now construct a set of $k$ generators $l_i = (T_i,S_i)$ for
$i=1,2, \ldots ,k$ using a total of $n$ carets
which we will show generate a subgroup of $F$
isomorphic to $H$.  We let $h_i = (T'_i,S'_i)$ as a tree pair diagram,
and $s = \sum_{i=1}^m |h_i|$. We let $(A,B)$ be a reduced pair of trees
with $n-(s+k)$ carets in each tree. We take  $n$ to be larger
than $s+k$ in order
 to construct $(A,B)$ in this way. We define $l_1$ by taking $T_1$ to
be the tree with a root caret whose left subtree is $T'_1$ and whose
right subtree is $A$.  Similarly, we let $S_1$ be the tree with a root
caret whose left subtree is $S'_1$ and whose right subtree is $B$.

For $2 \leq i \leq m$, we let $T_i$ be the tree consisting of a root
caret whose left subtree is $T'_i$ and whose right subtree is empty.
We let $S_i$ be the tree consisting of a root caret whose left subtree
is $S'_i$ and whose right subtree is empty.  For $m+1 \leq i \leq
k$, we let $l_i$ be the identity represented by a pair of trees each
containing a single caret.

We note that by construction, all tree pair diagrams constructed in
this way are reduced. We have $k$ root carets (counting one caret per pair), to which we attached $s$ carets for all the $(T_i',S_i')$ pairs, $n-(s+k)$ carets for the $(A,B)$ pair. This totals to $k+s+n-(s+k)
=n$ ensuring that the $k$-tuple constructed lies in the desired sphere.

It is clear that $\langle l_1,l_2, \ldots ,l_k \rangle$ generate a
subgroup of $H \times \Z$, where the isomorphic copy of $H$ lies in
the first factor of the standard $F \times F$ subgroup of $F$ and
where we take $(A,B)$ to be the generator of the $\Z$ factor which
lies in the second factor of the standard $F \times F$ subgroup.  We
now claim that $\langle l_1,l_2, \ldots ,l_k \rangle \cong H$.  We
use the coordinates $(h,t^a)$ on $H \times \Z$, where $h \in H$ and
$t = (A,B)$.  We define a homomorphism from $H \times \Z$ to $H$ by
taking the first coordinate of $(h,t^a)$.  When restricted to
$\langle l_1,l_2, \ldots ,l_k \rangle$, this map is onto by
construction.

To show this projection map is injective, we suppose that $(1, t^a)$
lies in the kernel, for $a \neq 0$.  Thus $\langle l_1,l_2, \ldots
,l_k \rangle \subset H \times \Z$ has a relator $\rho$ which, when
projected to $H$, yields a relator $r$ of $H$, and when considered
as a word in $\langle l_1,l_2, \ldots ,l_k \rangle$, has a second
coordinate not equal to the identity. But any relator $r$ of $H$,
when each $h_1$ is written as a word in $x_0$ and $x_1$, satisfies
$E_0(r) = 0$. Since the only generator of $H$ with $E_0(h_i) \neq 0$
is $h_1$, we see that $r$ must have the same number of $h_1$ and
$h_1^{-1}$ terms in it.  Thus $\rho$ must have the same number of
$l_1=(h_1,t)$ and $l_1^{-1}$ terms.  Since $l_1$ is the only
generator of $\langle l_1,l_2, \ldots ,l_k \rangle$ which can change
the $\Z$ coordinate of a product, having equal numbers of $l_1$ and
$l_1^{-1}$ terms in our relator $\rho$ implies that when the $H$
coordinate is the identity, the second coordinate must be $t^0$.
Thus projection to the first factor is an isomorphism when
restricted to $\langle l_1,l_2, \ldots ,l_k \rangle$, and we
conclude that this group is isomorphic to $H$.

We now show that the set of $k$-tuples of tree pair diagrams
constructed in this way is visible in \specs$(F)$.  There are
$r_{n-(s+k)}$ ways to choose the pair $(A,B)$, which had $n-(s+k)$
carets, and which determined the $l_1$ generator in this
construction. Thus we see that
$$
\lim_{n \rightarrow \infty} \frac{r_{n-(s+k)}}{
|\mathrm{Sph}^{\mathrm{sum}}_k(n)|} \geq \lim_{n \rightarrow \infty}
\frac{r_{n-(s+k)}}{ r_{n+k-1}} =\mu^{-(s-1+2k)}>0
$$
using Lemmas  \ref{lem:SizeSum} and \ref{lem:quotients}.
\ep

The probabilistic motivation for the definition of a visible
subgroup $H$ is that a set of $k$ randomly selected reduced pairs
of trees will generate a subgroup isomorphic to $H$ with nonzero
probability.  In the preceding proof, we were able to show that any
given $m$-generator subgroup is visible in \specs$(F)$ using a
$k$-tuple of pairs of trees consisting of one ``large" tree pair
diagram, $m-1$ ``small'' tree pair diagrams, and finally $k-m$
``tiny" tree pair diagrams representing the identity.

Given a subgroup $H$ of $F$, the estimate given above on a lower bound for the density of the isomorphism class of $H$ is small but positive.  It follows from the proof of Theorem \ref{thm:specs} that we obtain larger estimates of this lower bound when the original subgroup $H$ is generated by elements with small tree pair diagrams.   For example, the asymptotic density of the isomorphism
class of the subgroup $\Z$ is at least $\mu^{-5} \approx
\frac1{750000}$ in the set of all $2$-generator subgroups, since
$k=2$ and $\Z$ can be generated by $x_0$ which has size 2. For other
nontrivial subgroups, the construction in this proof will require
more carets and the lower bounds we obtain will be even smaller, but
always positive.

\sect{Subgroup spectrum with respect to the max
stratification}\label{sec:max}

We now begin to compute the subgroup spectrum with respect to a
different stratification, the ``max'' stratification, of the set of
all $k$-generator subgroups of $F$.  We again let $X_k$ be the set
of unordered $k$-tuples of reduced pairs of trees, and define the
sphere of size $n$ to be the collection of $k$-tuples in which the
maximum size of any component is $n$:
\[\mathrm{Sph}_k^{\mathrm{max}}(n)=\left\{(t_1,\ldots, t_k) \; | \; \max_{i \in \{1,2, \ldots ,k\}} \{|t_i|\}=n\right\}\]
For example,  the triple of tree pairs in Figure~\ref{fig:Z3}
 (once $h_1$ is reduced) lies in $\mathrm{Sph}_3^{\mathrm{max}}(4)$.
Defining spheres in this way induces the desired stratification of
$X_k$.

We define the density of a subset $T\subseteq X_n$ with respect to
 the max stratification by
\[\lim_{n\ra\infty} \frac{|T\cap \mathrm{Sph}_k^{\mathrm{max}}(n)|}{|\mathrm{Sph}_k^{\mathrm{max}}(n)|}\]
and  \specm to be the set of visible isomorphism classes of subgroups
of $F$ with respect to the max stratification. As noted at the
end of  Section~\ref{sec:sum}, the sum stratification   is biased
towards $k$-tuples of tree pair diagrams which contain multiple
copies of the identity and other ``small" pairs of trees having few
carets. Using the maximum number of carets in a tree pair diagram to
determine size seems to yield a more natural stratification.

We find strikingly different results when we compute \specm as
compared to \specs(F).  For example, we show that $\Z$ lies in
$\mathrm{Spec}^{\mathrm{max}}_1(F)$ but not in \specm for larger
values of $k$.

As in Section~\ref{sec:sum}, we must first obtain bounds on the size
of the sphere of radius $n$ with respect to the max stratification.  We will use
these bounds in the proofs below.  We begin with a lemma
about sums of $r_n$.

\begin{lem}[Sums of  $r_n$]
\label{lem:B}
For $n\geq 2$,
$\sum_{i=1}^{n-1} r_i \leq r_n$.
\end{lem}
\bp
Since $r_1=1<r_2=2$ the statement holds for $n=2$.
We assume for induction the statement is true for $k\geq 2$. Then
\begin{align*}
\sum_{i=1}^{k} r_i &= \sum_{i=1}^{k-1} r_i + r_k \leq 2r_k
\end{align*} by inductive assumption. We consider the set of  reduced tree pairs
with $k+1$ carets in each tree, where either the right child of each
root is empty, or the left child of each root is empty. In each case
there are $r_k$ ways to arrange the $k$ carets on the nonempty leaf,
and these tree pairs form disjoint subsets of the set of all reduced
pairs of trees with $k+1$ carets. Thus $2r_k\leq r_{k+1}$ which
completes the proof.
\ep

\begin{lem}[Size of $\mathrm{Sph}^{\mathrm{max}}_k(n)$]\label{lem:SizeMax}
For $k\geq 1$ and $ n\geq k$,
\[\frac1{k!}(r_n)^k\leq |\mathrm{Sph}^{\mathrm{max}}_k(n)| \leq k(r_n)^k\]
\end{lem}
\bp
For the lower bound, there are $(r_n)^k$ ordered $k$-tuples of
reduced
 tree pairs where {\em every} pair has $n$ carets. Since
 Sph$_k^{\mathrm{max}}(n)$ consists of unordered
 tuples then dividing this by $k!$ gives a lower bound.

For the upper bound, at least one of the $k$ tree pairs must have
$n$ carets. For $1\leq i\leq k$ suppose that $i$ tree pairs have
exactly $n$ carets, and the remaining $k-i$ tree pairs have strictly
less than $n$ carets. There are at most $(r_n)^{i}$
ordered
$i$-tuples of $n$-caret tree pairs, and so at most this many unordered $i$-tuples, and
 at most $\left( \sum_{j=1}^{n-1} r_j
\right)^{k-i}$
ordered
$(k-i)$-tuples of   tree pairs
with at most $n-1$ carets each,  and so at most this many unordered $(k-i)$-tuples.

So for each $i$ the number of unordered $k$-tuples of tree pairs
where $i$ pairs have $n$ carets and $k-i$ pairs have less than $n$
carets is at most
 \[{(r_n)^{i}}{\left( \sum_{j=1}^{n-1} r_j
\right)^{k-i}} \leq {(r_n)^{i}}{( r_n
)^{k-i}} ={(r_n)^k}\] by Lemma \ref{lem:B}.
Since our $k$-tuples of tree pairs are unordered, without loss of
generality we can list the ones containing $n$ carets first.

Thus for the total number of $k$-tuples,l we have at most
\[\sum_{i=1}^{k}(r_n)^k  =
k(r_n)^k.\]
\ep

We begin by showing that $\Z^k$ is present in \specm for all $k \geq 1$. We prove
that $\Z \notin$ \specm for $k > 1$, and conjecture that $\Z^m$ is
not visible in \specm for $k >m$.  In the proof below, we construct
a particular collection of subgroups of $F$ isomorphic to $\Z^k$,
all of whose generators have a common form, and show that this
collection of subgroups is visible.  Presumably, the actual density
of the isomorphism class of subgroups of $F$ isomorphic to $\Z^k$ is
considerably larger.

\begin{lem}[\specm is nonempty]
\label{lem:Z}
 $\Z^k \in$ \specm for all $k \geq 1$.
\end{lem}

\bp
We let $T$ be the tree consisting of a string of $k-1$ right carets. We
construct a set of $k$ pairs of trees which generate a subgroup of
$F$ isomorphic to $\Z^k$ as described in Section~\ref{sec:subgroups}.

We let $(A_i,B_i)$ be a reduced pair of trees each with $n-(k-1)$
carets for $i=1,2, \ldots ,k$. We let $h_i$ be the pair of trees
obtained by taking the pair $(T,T)$ and attaching $A_i$ to the
$i$-th leaf of the first copy of $T$, and $B_i$ to the $i$-th leaf
of the second copy of $T$.  We reduce the tree pair generated in
this way (which will be necessary for $i=1, \ldots, k-2$) to obtain
the reduced representative for $h_i$, which we again denote $h_i$.
We note that $h_k$ will have $n$ carets in each tree in its pair, so
this tuple does lie in the proper sphere of the stratification.
 As discussed above, the set
$\{h_1,h_2, \ldots ,h_k\}$ will generate a subgroup of $F$
isomorphic to $\Z^k$.

We compute the  density of the set of $k$-tuples of pairs
of trees constructed in this way to be at least:
$$\lim_{n \rightarrow \infty} \frac{(r_{n-k+1})^k}{  k (r_n)^k} =\frac{1}{k}\mu^{-k^2+k}>
0$$ using Lemma \ref{lem:quotients} and the upper bound from
Lemma~\ref{lem:SizeMax}.  Thus $\Z^k$ is visible in \specm.
\ep

For example, this shows that the density of $\Z^2$ in the set of 2-generator
subgroups is at least $\frac12 \mu^{-2} \approx \frac1{500}$.

We now show that a subgroup $H$ of $F$ cannot appear in \specm for
values of $k$ smaller than the rank of the abelianization $H_{ab}$.
\begin{lem}[Abelianization]
\label{lem:abelianization} We let $H$ be a subgroup of $F$, and let
$n$ be the rank of the abelianization $H_{ab}$ of $H$.  Then $H
\notin $ \specm for $k < n$.
\end{lem}
\bp
Since the rank of $H_{ab}$ is $n$, we know that $H$ cannot be
generated with fewer than $n$ elements.  Thus $H$ cannot be visible
in \specm for $k < n$.
\ep

Aside from straightforward obstructions like the group rank and the
rank of the abelianization, it is not clear what determines the
presence of an isomorphism class of subgroup in a given spectrum. In
general, it is difficult to show that an isomorphism class of
subgroup is not present in a particular spectrum.  This is because it can be
difficult to systematically describe all possible ways of generating
a subgroup isomorphic to a given one.  However, in the case of $\Z$,
we can show that $\Z$ is not present in the $k$-spectrum for $k\geq
2$. This highlights a major difference between the composition of
\specs$(F)$ and \specm, since $\Z$ appears in all spectra with
respect to the sum stratification. As a subgroup of $F$ with a
single generator is either the identity or infinite cyclic, it
follows that $\mathrm{Spec}^{\mathrm{max}}_1(F)$ contains only $\Z$.

\subsection{Proving that sum and max are not the same}
The goal of this section is to prove the following theorem.

\begin{thm}[$\Z$ not visible]\label{thm:Znotvisible}
With respect to the max stratification, the spectrum
$\mathrm{Spec}^{\mathrm{max}}_1(F) = \{\Z\}$ and for any $k \geq 2$,
we have that $\Z \notin$ \specm.
\end{thm}

The essence of this proof is that if $k$ group elements generate a subgroup isomorphic to $\Z$, then they must all be powers of a common element.  Thus we make precise the notion that counting the number of $k$-tuples which generate a subgroup isomorphic to $\Z$ is, up to a polynomial factor, the same problem as choosing a single reduced tree pair as the generator of the subgroup.

We begin with some elementary lemmas relating the slope of the first non-identity linear piece of an element $f \in F$ and the number of carets in the reduced tree pair diagram representing that element.

\begin{lemma}\label{lemma:carets}
If $f \in F$  has a break point with coordinates $(\frac{x}{2^m},\frac{y}{2^r})$ where $x,y$ are odd integers,  
 then the reduced tree pair diagram for $f$ has at least $\max(m,r)$ carets in each tree.
\end{lemma}

\bp
In each tree in the tree pair diagram, carets at level $k$ correspond to points in $[0,1]$ with denominator $2^k$.  The lemma follows.
\ep

\begin{lemma}\label{lemma:slope}
Suppose  that the first non-identity linear piece of $f \in F$ has slope $2^r$ for $r \neq 0$.  Then the reduced tree pair diagram for $f$ has at least $|r|$ carets.
\end{lemma}

\bp
Suppose that the first non-identity linear piece of $f$ with slope $2^r$ has endpoints with coordinates $(\frac{a}{2^s},\frac{a}{2^s})$ and $(\frac{b}{2^x},\frac{c}{2^y})$ where $a,b,c$ are odd integers. 
We easily see that
$$2^r = \frac{\frac{c}{2^y}-\frac{a}{2^s}}{\frac{b}{2^x}-\frac{a}{2^s}}$$ Factoring out the highest power of $2$ possible from the denominator and the numerator of this fraction, and letting $m_1 = \min(s,x)$ and $m_2 = \min(s,y)$, we obtain
$$\frac{\frac{1}{2^{m_1}}(\frac{c}{2^{y-m_1}}-\frac{a}{2^{s-m_1}})}{\frac{1}{2^{m_2}}(\frac{c}{2^{y-m_2}}-\frac{a}{2^{s-m_2}})} = \frac{2^{m_2}}{2^{m_1}}A = 2^r$$
where no additional powers of $2$ can be factored out of the $A$ part of this expression.  Thus we see that one of $m_1, m_2 $ must be at least $|r|$, and thus it follows from Lemma \ref{lemma:carets} that the tree pair diagram for $f$ has at least $|r|$ carets.
\ep

We will use the coordinates of the first breakpoint to vastly over count the number of pairs of tree pair diagrams that we are considering.  However, even this vast over counting will work for the final argument.  We also need the following elementary lemma that follows from Lemma \ref{lemma:slope}.

\begin{lemma}\label{lemma:roots}
Let $f \in F$ have a reduced tree pair diagram with $n$ carets.  Then $f$ does not have an $m$-th root for $m>n$.
\end{lemma}

\bp
Suppose that $f$ has an $m$-th root $h$ for some $m>n$.   If $f$ is the identity on $[0, \epsilon]$, then any root or power of $f$ will be the identity on this interval as well.  Let the slope of the first non-identity linear piece of $h$ be $2^r$ for $r \neq 0$, and have left endpoint $\alpha = (\frac{a}{2^s},\frac{a}{2^s})$ for $a$ odd.
 Then the slope of $h^m$ near $\alpha$ is $2^{rm}$ and $|rm|>n$ since $m>n$.  Thus it follows from Lemma \ref{lemma:slope} that the tree pair diagram for $f = h^m$ has more than $n$ carets, a contradiction.
\ep

The proof of Theorem \ref{thm:Znotvisible} is divided into the following three lemmas.  Note that Lemma \ref{lemma:Z-spec2} is a special case of Lemma \ref{lemma:Z-speck}, but is included to illustrate the ideas involved.

\begin{lemma}\label{lemma:Z}
With respect to the max stratification, the spectrum
$\mathrm{Spec}^{\mathrm{max}}_1(F) = \{\Z\}$.
\end{lemma}
\bp
It follows from Lemma~\ref{lem:Z} that $\Z \in
 \mathrm{Spec}^{\mathrm{max}}_1(F)$. The only other possible
candidate for a subgroup isomorphism class in
$\mathrm{Spec}^{\mathrm{max}}_1(F)$ is that of the identity, and the
only reduced tree pair diagram representing the identity is of size
1. The number of reduced tree pairs representing the identity is
0 for size $n > 1$, and thus the density of the isomorphism class of
the identity subgroup when $k=1$ is $0$. We conclude that
$\mathrm{Spec}^{\mathrm{max}}_1(F) = \{\Z\}$.
\ep

To see that $\Z \notin$ \specm for any $k \geq 2$, we begin by
over counting the number of $k$-tuples of elements which can generate a subgroup isomorphic to $\Z$.

\begin{lemma}\label{lemma:Z-spec2}
For a fixed $n>1$, there are at most $(2n+1)(n+1) r_n$ distinct unordered pairs of elements $f,g \in F$ so that
\begin{enumerate}
\item the number of carets in each tree pair diagram is at most $n$,
\item the number of carets in at least one tree pair diagram is equal to $n$, and
\item $\langle f,g \rangle \cong \Z$.
\end{enumerate}
\end{lemma}

\bp
Since $\langle f,g \rangle \cong \Z$ we know that $f$ and $g$ are powers of a common element.  Note that this includes the case where this common element is either $f$ or $g$.  By assumption, one of $f$ and $g$ has $n$ carets in its tree pair diagram; without loss of generality we assume that it is $f$.  Thus there are $r_n$ choices for $f$.

From Lemma \ref{lemma:roots} we know that $f$ may have $m$-th roots for $0 \leq m \leq n$.  It follows from \cite{brin-squier}, Theorem 4.15 that if $f \in F$ has an $m$-th root, then that root is unique.  Denote the possible roots of $f$ by $r_0,r_1, \cdots ,r_j$ for $0 \leq j \leq n$.  Note that we are including $f$ itself as the $0$-th root.  We also know that $g$ must be a power of one of those (at most) $n+1$ possible roots, so there is an $i$ so that $g = r_i^a$ for some integer $a$.  Since $g$ has at most $n$ carets in its tree pair diagram, it follows from Lemma \ref{lemma:slope} that this exponent $a$ is at most $n$ in absolute value.  To see this, let $(x,y)$ be the first break point of $f$ so that the slope of the linear piece following $(x,y)$ is $2^{\alpha}$ for $\alpha \neq 0$.  Then it is easy to see that the slope to the right of $(x,y)$ in $f^k$ is $2^{k \alpha}$ and the statement then follows from Lemma \ref{lemma:slope}.  Thus there are $2n+1$ choices for the exponent $a$ so that $h^a = g$ since $|a| \leq n$.  In total, the number of ways we can construct a pair of this form is $(2n+1)(n+1)r_n$. Again, this count includes many pairs of elements that do not satisfy the requirements of the proposition, but all elements that do satisfy those conditions are counted in this argument.
\ep


\begin{lemma}\label{lemma:Z-speck}
For a fixed $n>1$, there are at most $(2n+1)^{k-1}(n+1) r_n$ distinct unordered $k$-tuples of elements $f_1,f_2, \cdots ,f_k  \in F$ so that
\begin{enumerate}
\item the number of carets in each tree pair diagram is at most $n$,
\item the number of carets in at least one tree pair diagram is equal to $n$, and
\item $\langle f_1,f_2, \cdots ,f_k \rangle \cong \Z$.
\end{enumerate}
\end{lemma}

\bp
The argument follows the proof of Lemma \ref{lemma:Z-spec2}.  There must be some element $h$ which generates this copy of $\Z$, that is, all $f_i$ are powers of this element $h$.  Suppose without loss of generality that $f_1$ has $n$ carets in its tree pair diagram.  Then $f_1$ may have $m$-th roots for $0 \leq m \leq n$, which we denote  $r_0,r_1, \cdots ,r_j$ for $0 \leq j \leq n$.  The same reasoning shows that for each $i$ we must have $f_i = r_j^{e_j}$, where $0 \leq j \leq n$ and $|e_j| \leq n$, where the latter inequality follows from Lemma \ref{lemma:slope}.  We then see that the number of such $k$-tuples is $(2n+1)^{k-1}(n+1) r_n$.
\ep

We now finish the proof of Theorem \ref{thm:Znotvisible}.

\noindent
{\em Proof of Theorem \ref{thm:Znotvisible}.}  For $k \geq 2$ we see that the density of $k$-tuples of pairs of trees which generate a subgroup isomorphic to $\Z$ is 
$$\lim_{n \rightarrow \infty} \frac{(2n+1)^{k+1} (n+1) r_n}{\frac{1}{k!} (r_n)^k} = 0$$
using the bound on the size of the $n$-sphere in the max stratification given in Lemma \ref{lem:SizeMax} as well as the upper bounds proven in Lemmas \ref{lemma:Z-spec2} and \ref{lemma:Z-speck}. The first statement in the theorem follows from Lemma \ref{lemma:Z} and the second from the above limit.
\ep

We note that this approach does not appear to generalize to show that $\Z^m$ is
not visible in \specm for $k>m$, as it is difficult to recognize when a collection
of tree pair diagrams generates a subgroup isomorphic to $\Z^m$ for $m\geq2$.

\subsection{Further results within the max stratification}
Apart from $\mathrm{Spec}^{\mathrm{max}}_1(F)$, it 
seems quite difficult to
compute the complete list of subgroups 
which appear
in \specm. Indeed, ignoring any consideration of densities, a
complete list of even the $2$-generated subgroups of $F$ is not
known (see \cite{\ProbList} Problem 2.4). For $k=2$ we can say the
following.
\begin{prop}[2-spectrum of $F$]
\label{prop:spec2} Let $H = \langle h_1,h_2 \rangle$ be a subgroup
of $F$.  Then either $H$ or $H \times \Z$ lies in
$\mathrm{Spec}^{\mathrm{max}}_2(F)$.  If $H_{ab} \cong \Z \bigoplus
\Z$, then $H \in \mathrm{Spec}^{\mathrm{max}}_2(F)$, otherwise $H
\times \Z \in \mathrm{Spec}^{\mathrm{max}}_2(F)$.
\end{prop}

\bp We may assume, quoting Lemmas~\ref{lem:outofcommutator}
and \ref{lem:exponentsum} that
if $H = \langle h_1,h_2 \rangle$ that
\bi
\item $h_1 \notin [F,F]$
\item when $h_1$ is expressed as a word in $x_0$ and $x_1$, the
exponent sum of all the instances of $x_0$ is not equal to $0$, and
\item when $h_2$ is expressed as a word in $x_0$ and $x_1$, the
exponent sum of all the instances of $x_0$ is equal to $0$.
\ei
As tree pair diagrams, we use the notation $h_i = (S_i,T_i)$.

\begin{figure}[ht!]
  \bc
\includegraphics[ scale=.5]{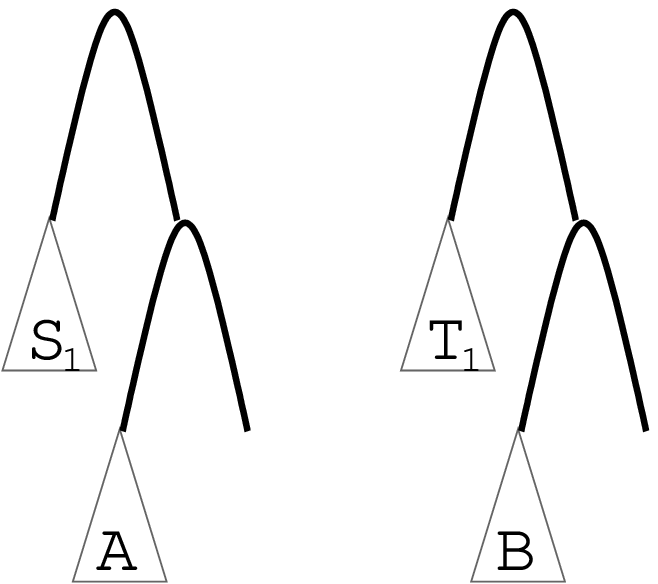}
\hspace{12mm}
   \includegraphics[ scale=.5]{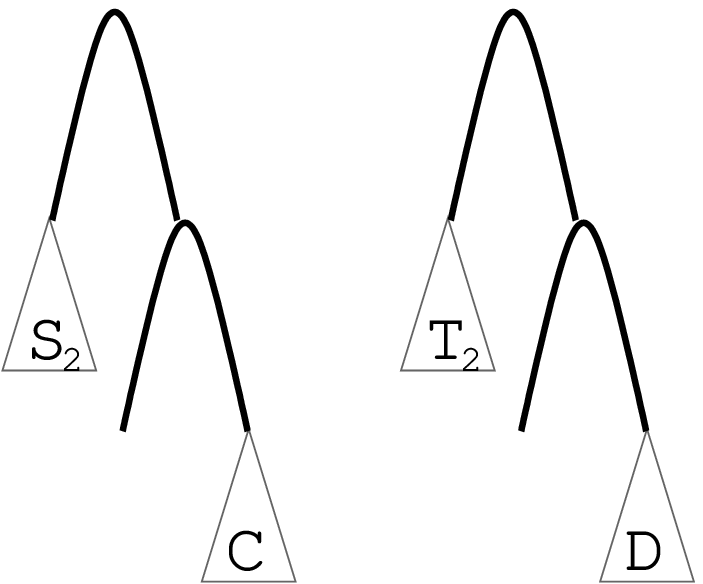}
\ec
  \caption{Constructing the tree pairs $k_1, k_2$ which generate a subgroup of $H\times\Z$.}
  \label{fig:2spec}
\end{figure}

We create a new set of generators $k_1=(X_1,Y_1)$ and
$k_2=(X_2,Y_2)$ for a two generator subgroup of $F$ as follows. We
let $T$ be the tree consisting entirely of two right carets, whose
leaves are numbered $1,2$ and $3$, and let $(A,B)$ and $(C,D)$  be
arbitrary reduced pairs of trees so that $(A,B)$ has $n-N(h_1)-2$
carets in each tree and $(C,D)$ has $n-N(h_2)-2$ carets in each
tree. We construct $X_1$ by attaching $S_1$ to leaf $1$ of $T$ and
$A$ to leaf $2$ of $T$. We construct $Y_1$ by attaching $T_1$ to
leaf $1$ of $T$ and $B$ to leaf $2$ of $T$.  We construct $X_2$ by
attaching $S_2$ to leaf $1$ of $T$ and $C$ to leaf $3$ of $T$. We
construct $Y_2$ by attaching $T_2$ to leaf $1$ of $T$ and $D$ to
leaf $3$ of $T$, as in Figure~\ref{fig:2spec}. Note that each tree
has size $n$, and we assume without loss of generality that $n >
\max\{N(h_i)\}+4$ so that the trees $A,B,C$ and $D$ each have at least two carets.

One may easily verify that  $k_1$ and $k_2$ generate a subgroup of the
standard $F \times F$ subgroup in which the subgroup you obtain on
the first factor of $F$ is simply $H$.  Also, $t=(A,B)$ and
$s=(C,D)$ each generate a copy of $\Z$ in the second factor of $F
\times F$ provided that neither tree pair diagram represents the
identity. Let $K \cong \langle k_1,k_2 \rangle$.  Then by
construction, $K \subset H \times \Z^2$, where the first $\Z$ is
generated by $t=(A,B)$ and the second by $s=(C,D)$.

We first show that the set of subgroups $K$ constructed in this way
is visible in $\mathrm{Spec}^{\mathrm{max}}_2(F)$, and then we
discuss of what isomorphism class of subgroups we have constructed
using these elements. By Lemmas~\ref{lem:SizeMax} and
\ref{lem:quotients} the
 density of pairs of tree pair diagrams constructed in
this way is at least
$$\lim_{n \rightarrow \infty} \frac{\left(r_{n-N(h_1)-2}\right)\left(r_{n-N(h_2)-2}\right)}{2 (r_n)^2} =
 \frac12\mu^{-N(h_1)-N(h_2)-4}>
0.$$

We claim that $K$ is either isomorphic to $H$ or to $H \times \Z$.
Use the coordinates $(w,t^a,s^b)$ on $H \times \Z^2$ where $w \in
H$.  It is easy to see that for every element $h \in H$, there is at
least one $k \in K$ represented by the coordinates $(h,t^a,s^b)$ for
some $a,b \in \Z$.  We first show that for each $h \in H$, there is
a unique second coordinate. Suppose that $w_1=(h,t^a,s^b)$ and
$w_2=(h,t^c,s^d)$ both lie in $K$, and thus the product $w_1w_2^{-1}=(Id,
t^{a-c},s^{b-d})$ also lies in $K$. Thus there is some relation
$\rho$ in $H$ expressed in terms of $h_1$ and $h_2$ so that when we
replace $h_i$ with $k_i$ we obtain the element $(Id,
t^{a-c},s^{b-d}) \in K$.  Since the generator $t$ of $\Z$ is linked
to $h_1$ in $k_1$, and the $t$ coordinate of $(Id, t^{a-c},s^{b-d})$
is not zero, we conclude that in $\rho$, the exponent sum of all
instances of the generator $h_1$ is not equal to zero.

Recall that $h_1$ was chosen so that when $h_1$ is expressed as a
word in $x_0$ and $x_1$, the exponent sum of all the instances of
$x_0$ is not equal to $0$, but $h_2$ does not have this property.
Any relation in $H$ can be written in terms of $x_0$ and $x_1$ to
yield a relation of $F$, and thus any relation in $H$ must have the
total exponent sum of all instances of $x_0$ equal to $0$. By our
choice of $h_1$ and $h_2$, we see that a relation of $H$ {\em must}
have the exponent sum of all instances of the generator $h_1$ equal
to zero.  Thus we must have $a=c$ in our coordinates above.

We have now shown that either $K \cong H \times \Z$ or $K \cong H$.
Suppose that $H_{ab} \cong \Z \bigoplus \Z$.  Then $(H \times
\Z)_{ab} \cong \Z^3$ and it follows from
Lemma~\ref{lem:abelianization}
 that $H \times \Z \notin \mathrm{Spec}^{\mathrm{max}}_2(F)$.  In this
case we must have $K \cong H$.

Suppose that $H_{ab} = \Z$.  In this case, either $h_2 \in [H,H]$ or
$h_1^c = h_2^dw$ for some non-identity element $w \in [H,H]$ and
integers $c,d$. In either case, there is a relator of $H$ in which
the total exponent sum on the instances of $h_2$ is nonzero. Since
the $s$ coordinate of the second $\Z$ factor in $H \times \Z^2$ is
linked to the $h_2$ generator in $k_2$, there is a way to realize
both $(h,t^a,s^b)$ and $(h,t^a,s^d)$ in $K$ with $d \neq b$. Thus we
must have $K \cong H \times \Z$.
\ep

It follows from Proposition \ref{prop:spec2} that
$\mathrm{Spec}^{\mathrm{max}}_2(F)$ contains $\Z^2, F,$ and
$\Z\wr\Z$, and from Theorem \ref{thm:Znotvisible} that it does not
contain $\Z$ or $\{Id\}$.

We have seen above that it can be difficult to ascertain when a
particular isomorphism class of subgroup is present in a given
spectrum. Furthermore, the example of $\Z$ shows that presence in a
given spectrum does not necessarily imply presence in spectra of
higher index.

We find that $F$ is a very special two generator subgroup of itself,
and exhibits behavior unlike that of $\Z$. As long as $k
\geq 2$, we can show that $F \in $ \specm.  We call this behavior
{\em persistence}; that is, a subgroup $H$ is {\em persistent} if there is
an $l$ so that $H \in $ \specm for all $k \geq l$. In the small set
of groups whose spectra have been previously studied, no subgroups
have shown this persistent behavior.  As noted in the introduction, the current
known examples of subgroup
spectra all find that the free group $F_k$ is generic in the
$k$-spectrum. In Thompson's group $F$, we find a wealth of examples
of this persistent behavior. In the previous section, we effectively proved that {\em every} non-trivial finitely generated subgroup of $F$ is persistent with respect to the sum stratification (Theorem \ref{thm:specs}). As a corollary of
Theorem~\ref{thm:Fpersis} below and the techniques in
 Lemma~\ref{lem:Z} above, it will follow that $F^n \times \Z^m$ and  $F^n
\wr \Z^m$ are also persistent with respect to the {\em max} stratification, with $l=2n+m$.

\begin{thm}[F is persistent]\label{thm:Fpersis}
$F$ lies in \specm for all $k \geq 2$.
\end{thm}

\bp
Since $F$ can be generated by two elements, and $F_{ab} \cong \Z
\bigoplus \Z$, it follows from Proposition~\ref{prop:spec2} that $F
\in \mathrm{Spec}^{\mathrm{max}}_2(F)$.  We now show that $F \in$
\specm for all
 $k> 2$.

\begin{figure}[ht!]
  \bc
\includegraphics[ scale=.45]{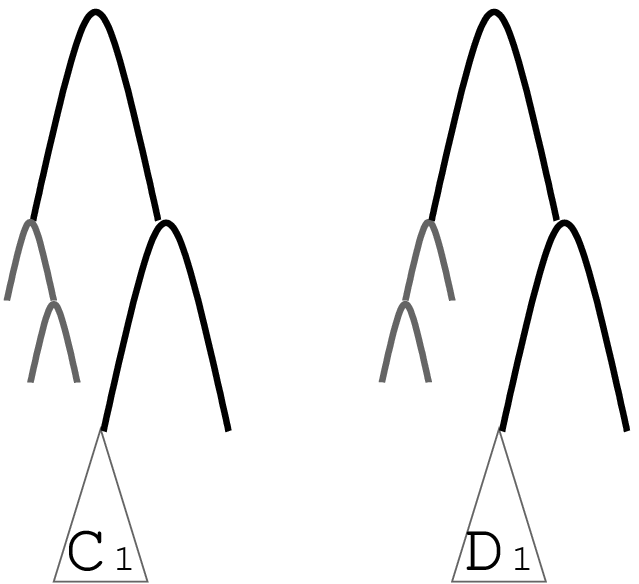}
\hspace{12mm}
   \includegraphics[ scale=.45]{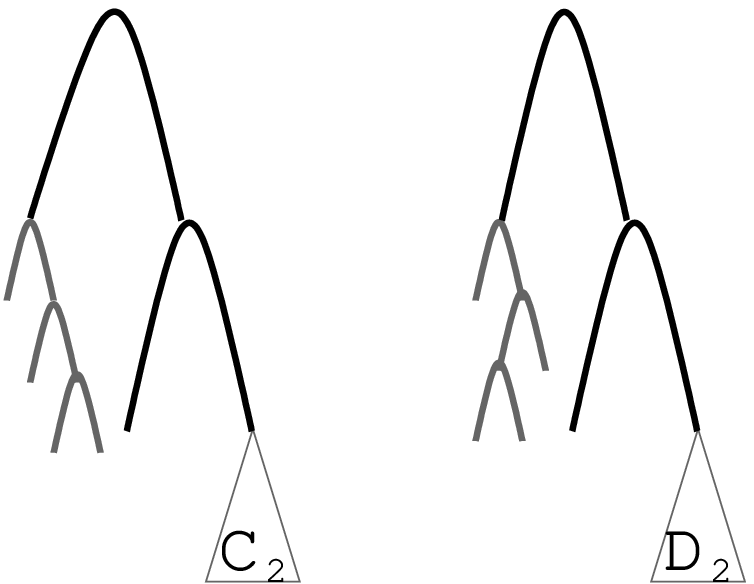}
  \hspace{12mm}
  \includegraphics[ scale=.45]{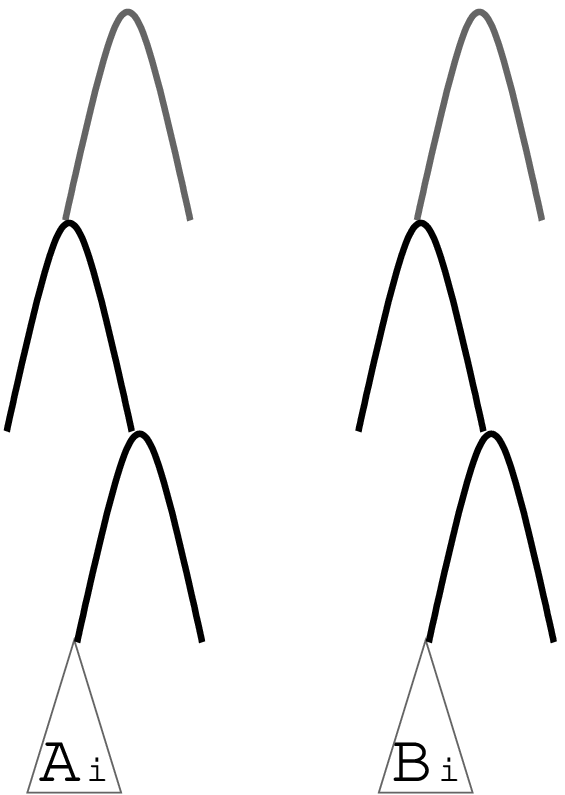}
\ec
  \caption{Constructing the tree pairs $h_1, h_2$, and $h_i$ generating a subgroup of $F\times\Z^2$.}
  \label{fig:persis}
\end{figure}

We define $k$ generators $h_1,h_2, \ldots ,h_k$ which generate a
subgroup of $F$ isomorphic to $F$, in such a way that the set of
$k$-tuples pairs of trees of this form is visible.  As reduced tree
pair diagrams, we use the notation $h_i = (T_i,S_i)$.  We begin by
defining $h_1$ and $h_2$. We let $x_0 = (T_{x_0},S_{x_0})$ and $x_1
= (T_{x_1},S_{x_1})$ as tree pair diagrams,  $(C_1,D_1)$ any reduced
pair of trees with $n-4$ carets in each tree and $(C_2,D_2)$ any
reduced pair of trees with $n-5$ carets in each tree. We let $T$ be
the tree with two right carets, and three leaves numbered $1,2,3$.
We construct $h_1$ and $h_2$ as follows: \bi
\item We let $T_1$ be the
tree $T$ with $T_{x_0}$ attached to leaf $1$ and $C_1$ attached to
leaf $2$.
\item We let $S_1$ be the tree $T$ with $S_{x_0}$ attached to leaf $1$
and $D_1$ attached to leaf $2$.
\item We let $T_2$ be the tree $T$ with
$T_{x_1}$ attached to leaf $1$ and $C_2$ attached to leaf $3$.
\item We let $S_2$ be the tree $T$ with $S_{x_1}$ attached to leaf $1$
and $D_2$ attached to leaf $3$.
\ei
This construction is shown in  Figure~\ref{fig:persis}.

For fixed $n$,  let $(A_i,B_i)$ be any reduced pair of trees with
$n-3$ carets for $i=3,4, \ldots ,k$.
Note that there are $r_{n-3}$ ways to choose each such pair.
Construct a reduced $(n-1)$-caret tree pair that represents an element of $[F,F]$ by
attaching the pair $(A_i,B_i)$ to a 2-caret tree as in Figure~\ref{fig:commutator}
in Section~\ref{sec:commutator}.
Call this pair $(A'_i,B'_i)$.
   We now define
$h_i = (T_i,S_i)$ for $i=3,4, \ldots ,k$ as follows:
 \bi
\item let $T_i$ consist of a root caret with $A'_i$ attached to its left
leaf, and
\item let $S_i$ consist of a root caret with $B'_i$ attached to its left
leaf.
\ei

The subgroup generated by the $\{h_i\}$ is clearly a subgroup of $F
\times \Z^2$, since the subtrees of the $h_i$ which are the left
children of the root carets, when taken as independent tree pair
diagrams, clearly generate a subgroup $H$ which is isomorphic to
$F$, as they contain the tree pair diagrams for $x_0$ and $x_1$.

Any relator which is introduced into $H$ by the inclusion of the
commutators $(A'_i,B'_i)$ as generators must hold true in $F$ as
well.  Since all relators of $F$ are commutators or conjugates of
commutators, all relators have exponent sum on all instances of
either $x_0$ and $x_1$ equal to zero.  Additionally, we know that
$x_0$ and $x_1$ are not commutators themselves.  Thus any new
relators introduced into $H$ by the inclusion of the commutators
$(A'_i,B'_i)$ as generators must also have  exponent sum on all
instances of either $x_0$ and $x_1$ equal to zero.  Using the
coordinates $(w,t^a,s^b)$ for elements of $H$, where $w \in F$, $t =
(C_1,D_1)$ and $s=(C_2,D_2)$, the argument given in
Proposition~\ref{prop:spec2} goes through exactly to show that $w
\in F$ has unique second and third coordinates, and thus $H \cong
F$.

To see that the set of $k$-tuples constructed in this way is
visible, note that the number of ways to construct them is
$r_{n-4}r_{n-5}(r_{n-3})^{k-2}$.  The choices are in the $C_i,D_i$ trees  which generate
$\Z^2$, and the $A_i',B_i'$ trees which are used to construct elements of
$[F,F]$. Thus we compute the  density of this set of
$k$-tuples to be at least
$$\lim_{n \rightarrow \infty} \frac{r_{n-4}r_{n-5}(r_{n-3})^{k-2}}{k (r_n)^k}
= \frac1{k}\mu^{-4}\mu^{-5}(\mu^{-3})^{k-2}
>
0$$ by Lemmas~\ref{lem:SizeMax} and \ref{lem:quotients}.
\ep

This proof used two very special properties of the whole group $F$
which are not generally true for subgroups of $F$.  First, there is
an explicit way of characterizing tree pair diagrams corresponding
to elements in the commutator subgroup $[F,F]$, which allows us to
construct commutators containing a large arbitrary tree.  Second,
the relators of $F$ are all commutators themselves, and thus
including additional commutators as generators yields relators with
the appropriate exponent sums on $x_0$ and $x_1$.  Thus we do not
expect this persistent behavior from many other subgroups
of $F$.  However, we can adapt the ideas used above to prove that if
a subgroup $H$ of $F$ is visible in a particular spectrum, \specm,
then both the product $H \times \Z$ and the wreath product $H \wr \Z$ are
visible in $\mathrm{Spec}^{\mathrm{max}}_{k+1}(F)$.  As a corollary
of this fact and Theorem~\ref{thm:Fpersis}, we find that subgroups
which contain $F$ as a factor are indeed persistent.  We first need
the following straightforward lemma about densities of visible subgroups.

\begin{lem}
\label{lem:tech}
We let $H_k(n)$ denote the
 set of all $k$-tuples of tree pair diagrams which
generate a subgroup of $F$ isomorphic to $H$ with a maximum of $n$
carets in any pair of trees, such that at least one coordinate
realizes this maximum.
 If a subgroup $H$ is visible in
\specm then
\[\lim_{n \rightarrow \infty} \frac{|H_k(n)|}{(r_n)^k} \geq \lambda_k\]
for some $\lambda_k \in (0,1]$.
\end{lem}
\bp
\begin{align*}
\lim_{n \rightarrow \infty} \frac{|H_k(n)|}{(r_n)^k}
&\geq \lim_{n \rightarrow \infty} \frac{|H_k(n)|}{k! |\mathrm{Sph}^{\mathrm{max}}_k(n)| }\\
\end{align*}
by Lemma~\ref{lem:SizeMax}. Since $H$ is visible this limit equals
the density of $H$ with respect to the max stratification, and is
positive, which gives the result.
\ep

\begin{prop}[Closure under products]
\label{prop:subgroups}
If $H\in$
 \specm then $H\times \Z$ and $H\wr\Z$ lie in
$\mathrm{Spec}^{\mathrm{max}}_{k+1}(F)$.
\end{prop}

\bp
We construct the $k+1$ generators necessary to obtain a family of
subgroups of $F$ isomorphic to $H \times \Z$ in such a way that the
set of $(k+1)$-tuples of this form is visible.  The techniques are
similar to those used above.

We let $h_1,h_2, \ldots h_k$ be a set of $k$ generators for $H$.  We
will construct a set $l_1,l_2, \ldots ,l_{k+1}$ of generators for $H
\times \Z$.  We let $h_i = (T_i',S_i')$ as a reduced pair of trees, and
we must define $l_i = (T_i,S_i)$.  For $i=1, \ldots ,k$ we let $T_i$
consist of a root caret with $T_i'$ as its left subtree, and $S_i$
consist of a root caret with $S_i'$ as its left subtree. We let $(A,B)$
be a reduced pair of trees with $n-1$ carets.  To define $l_{k+1}$,
let $T_{k+1}$ consist of a root caret with $A$ as its right subtree,
and $S_{k+1}$ consist of a root caret with $B$ as its right subtree.

It is clear that the set $\{l_i\}$ generate a subgroup of $F$
isomorphic to $H \times \Z$.  We now show that the set of
$(k+1)$-tuples constructed in this way is visible in
$\mathrm{Spec}^{\mathrm{max}}_{k+1}(F)$.

To compute the  density of
the set of $(k+1)$-tuples constructed in this way which generate a
subgroup of $F$ isomorphic to $H \times \Z$, we compute the
following limit.
\begin{align*}
 \lim_{n \rightarrow \infty} \frac{|H_k(n-1)|r_{n-1}}{|\mathrm{Sph}^{\mathrm{max}}_{k+1}(n)|} &
\geq \lim_{n \rightarrow \infty}\frac{|H_k(n-1)|r_{n-1}}{(k+1)(r_n)^{k+1}}\\
\intertext{by Lemma~\ref{lem:SizeMax}} &=\lim_{n \rightarrow
\infty}\frac1{k+1}
\frac{|H_k(n-1)|}{(r_n)^k}\frac{r_{n-1}}{r_n} \\
&=\lim_{n \rightarrow \infty}\frac1{k+1}
\frac{|H_k(n-1)|}{(r_{n-1})^k}\frac{(r_{n-1})^k}{(r_n)^k}\frac{r_{n-1}}{r_n} \\
&\geq \frac{ \lambda_k \mu^{-k-1}}{k+1}> 0. \intertext{by
Lemmas~\ref{lem:tech} and \ref{lem:quotients}.}
\end{align*}

To see that $H\wr\Z$ lies in $\mathrm{Spec}^{\mathrm{max}}_{k+1}(F)$
under the same assumption on $H$, we construct slightly different
generators, and make an argument analogous to that in
Theorem~\ref{thm:Fpersis}.  As above, we let $h_1,h_2, \ldots h_k$
be a set of $k$ generators for $H$. We will construct a set
$l_1,l_2, \ldots ,l_{k+1}$ of generators which will generate a
subgroup of $(H \wr \Z) \times \Z$ which we show to be isomorphic to
$H \wr \Z$.

Let $H_k(n-3)$ be the set of all $k$-tuples which generate a
subgroup of $F$ isomorphic to $H$, where at least one tree pair
contains $n-3$ carets.  Let $\{h_i = (T_i',S_i')\} \in H_k(n-3)$.
Since $H$ is visible in \specm, Lemma~\ref{lem:tech} implies that
$$\lim_{n \rightarrow \infty} \frac{|H_k(n-3)|}{r_{n-3}^k} > 0.$$

\begin{figure}[ht!]
  \bc
\includegraphics[ scale=.45]{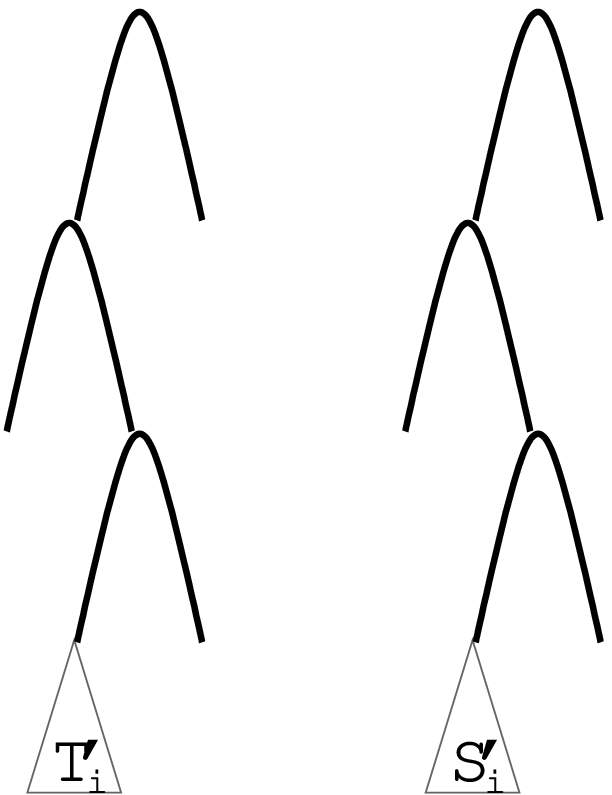}
\hspace{12mm}
   \includegraphics[ scale=.45]{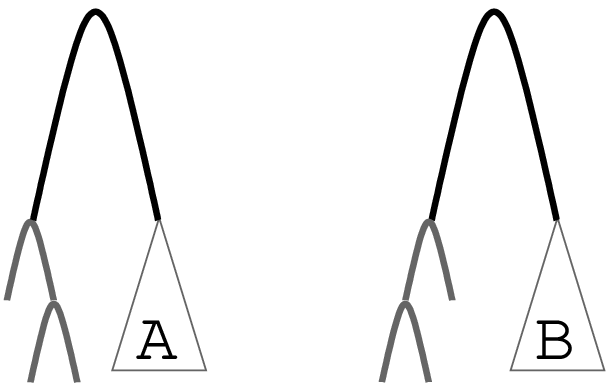}
\ec
  \caption{Constructing the pairs $l_i$  and $l_{k+1}$ generating a subgroup of $(H\wr\Z)\times\Z$.}
  \label{fig:wreath}
\end{figure}

We define $l_i = (T_i,S_i)$ for
$i=1,2, \ldots ,k+1$ as follows.
We let $T$ be the tree with two left carets, and one interior caret
attached to the right leaf of the caret which is not the root.
Number the leaves of $T$ by $1,2,3,4$.  For $i = 1,2, \ldots ,k$,
let $T_i$ be the tree $T$ with $T_i'$ attached to leaf $2$.  We let
$S_i$ be the tree $T$ with $S_i'$ attached to leaf $2$.  We let $(A,B)$
be any reduced pair of trees with $n-3$ carets.  We let $x_0 =
(T_{x_0},S_{x_0})$.  We define $l_{k+1}$ by taking $T_{k+1}$ to be a
single root caret with $T_{x_0}$ attached to its left leaf and $A$
attached to its right leaf.  We let $S_{k+1}$ be a single root caret
with $S_{x_0}$ attached to its left leaf and $B$ attached to its
right leaf.
See Figure~\ref{fig:wreath}.

 It is clear by the construction of our generators that
any element of $H \wr \Z$ can appear as the pair of left subtrees of
the root carets in any element of $\langle h_i \rangle$.  However,
we must show that $\langle h_i \rangle$ generates a subgroup of $F$
isomorphic to $H \wr \Z$ and not $(H \wr \Z) \times \Z$.  To do
this, we note that since $H \wr \Z$ is a wreath product, all
relators are commutators.  Thus the argument in Theorem~\ref{thm:Fpersis}
can be applied to show that $\langle h_i \rangle \cong H \wr \Z$
rather than $(H \wr \Z) \times \Z$.

We must now show that the set of $(k+1)$-tuples generated in this
way is visible in $\mathrm{Spec}^{\mathrm{max}}_{k+1}(F)$.  We let
$H_k(n)$ be the set of all $k$-tuples of tree pair diagrams which
generate a subgroup of $F$ isomorphic to $H$ with a maximum of $n$
carets in any pair of trees, such that at least one coordinate
realizes this maximum.
  The
 density of the set of $(k+1)$-tuples constructed in this
way which generate a subgroup of $F$ isomorphic to $H \wr \Z$ is
computed as follows. We have $r_{n-3}$ choices for the pair $(A,B)$,
and $|H(n-3)|$ is the number of $(T_i',S_i')$ generating sets for
$H$ with a maximum of $n-3$ carets in some pair. So together the
density is
\begin{align*}
\lim_{n \rightarrow
\infty}\frac{|H(n-3)|r_{n-3}}{|\mathrm{Sph}^{\mathrm{max}}_{k+1}(n)|}
&\geq  \lim_{n \rightarrow \infty}\frac1{k+1}\frac{|H(n-3)|r_{n-3}}{(r_n)^kr_n}\\
&\geq  \lim_{n \rightarrow \infty}\frac1{k+1}\frac{|H(n-3)|}{(r_{n-3})^k}\frac{(r_{n-3})^k}{(r_n)^k}\frac{r_{n-3}}{r_n}\\
&=\frac{1}{k+1}\lambda_k\mu^{-3k-3}>0
\end{align*}
{by Lemmas~\ref{lem:tech} and \ref{lem:quotients}.}

\ep

This proposition combined with Theorem~\ref{thm:Fpersis} allows us to find many isomorphism classes of subgroups in $\mathrm{Spec}^{\mathrm{max}}_k(F)$ for the appropriate value of $k$.








\bi\item The $l$-fold iterated wreath product of $\Z$ with
itself $\Z \wr \cdots \wr \Z$ lies in
$\mathrm{Spec}^{\mathrm{max}}_l(F)$.
\item
If $H$ is a persistent subgroup present in
 \specm for $k \geq l$, then $H\times \Z$ and $H\wr\Z$ are
 persistent subgroups present in
\specm for $k \geq l+1$.
\item
For $n \geq 1, m \geq 0$, and for all $k \geq
2n+m$, we have that $F^n \times \Z^m$ lies in \specm.

This shows that it is possible to have a
subgroup $H$ of $F$ so that both $H$ and $H \times \Z$ are contained
in the \specm for the same value of $k$; we can take $H = F^n \times
\Z^m$ and $k> 2m+n$.
\item
$F^n \wr \Z$ lies in \specm for $n\geq 1$ for all $k \geq 2n+1$.\ei

More generally, we can see that persistent subgroups can ``absorb''
 visible subgroups to form new persistent subgroups.

\begin{thm}[Products with persistent subgroups are
persistent]\label{thm:persistent-products} If $H$ is a subgroup
which is present in \specm and $K$ is a persistent subgroup which is
present in $\mathrm{Spec}^{\mathrm{max}}_l(F)$ for $l \geq l_0$,
then $H \times K$ is persistent and present in
$\mathrm{Spec}^{\mathrm{max}}_l(F)$ for $l \geq l_0+k$.
\end{thm}

\bp
Let $H_k(n)$ denote the set of all $k$-tuples of tree pair diagrams
which generate a subgroup of $F$ isomorphic to $H$ with a realized
maximum of $n$ carets in some coordinate.  Since $H \in$ \specm we
know from Lemma \ref{lem:tech} that $$\lim_{n \rightarrow \infty}
\frac{|H_k(n)|}{(r_n)^k} \geq \lambda_k$$ for some $\lambda_k \in
(0,1]$.

Let $K_l(n)$ denote the set of all $l$-tuples of tree pair diagrams
which generate a subgroup of $F$ isomorphic to $K$ with a realized
maximum of $n$ carets in some coordinate.  Since $K$ is persistent,
we know that for any $l \geq l_0$, the limit $$\lim_{n \rightarrow
\infty} \frac{|K_l(n)|}{(r_n)^l} \geq \lambda_l$$ for some
$\lambda_l \in (0,1]$.

Let $m=k+l$ for any $l \geq l_0$.  Form a generating set $\{t_1,t_2,
\cdots ,t_m\}$, where $t_i = (T_i,S_i)$, for $H \times K$ as
follows. Take any $k$-tuple $\delta \in H_k(n)$, where $\delta_i \in
\delta$ is represented by the pair of trees
$(T^{\delta}_i,S^{\delta}_i)$. Take any $l$-tuple $\eta \in K_l(n)$,
where $\eta_j \in \eta$ is represented by the pair of trees
$(T^{\eta}_j,S^{\eta}_j)$.
\begin{itemize}
\item For $1 \leq i \leq k$, let $T_i$ consist of a root caret with left
subtree $T_i^{\delta}$, and let $S_i$ consist of a root caret with left
subtree $S_i^{\delta}$.
\item For $k+1 \leq i \leq m$, let $T_i$ consist of a root caret with right
subtree $T_i^{\eta}$, and let $S_i$ consist of a root caret with right
subtree $S_i^{\eta}$.
\end{itemize}
This set of tree pairs generates a subgroup of $F
\times F$ isomorphic to $H \times K$. A lower bound on the density
of the isomorphism class of $H \times K$ is given by the following
positive valued limit:
$$\lim_{n \rightarrow \infty} \frac{|H_k(n)||K_l(n)|}{r_n^{k+l}} = \lim_{n \rightarrow \infty}
\frac{|H_k(n)|}{r_n^k}\frac{|K_l(n)|}{r_n^l} \geq \lambda_k
\lambda_l
>0.$$
~\ep

Thus, our analysis shows that the following subgroups are present in the $k$-spectrum
with respect to the max stratification:
\bi
\item The persistent subgroups $F$, $F \times F$, \ldots $F^n$ for $2n \leq
k$.
\item The persistent subgroups $F^n \times \Z^m$, for $2n+m \leq k, n \geq
1$.
\item The persistent subgroups $F^n \wr \Z$ for $2n+1 \leq k, n \geq 1.$
\item The abelian subgroup $\Z^k$ and the $k$-fold iterated product of $\Z$ with
itself.
\item The mixed direct and wreath products of $\Z$ with itself with $k$ terms,  including
for example $\Z^{k-1} \wr \Z$ and $(\Z \wr \Z \wr \Z) \times \Z^{k-3}$.
\item Various mixed direct and wreath products with $\Z$ such as
 $(F^2 \times \Z^3) \wr \Z \times \Z$ which is present in all $k \geq 9$, for example.

\ei

While the isomorphism classes of subgroups described above occur with positive densities in \specm for appropriate $k$, the lower bounds on their densities are very small.  In fact, the lower bound on the sum of the densities of all of these isomorphism classes of subgroups amounts to much less than $1\%$ of all isomorphism classes of subgroups in \specm.

We conclude with an open question about the isomorphism type of a random subgroup of the other Thompson's groups $T$ and $V$.  Although these groups contain $F$ as a proper subgroup, unlike $F$ they also contain free subgroups or rank $2$ and above.  What is the density of the set of free subgroups of a given rank within \spec$_k(T)$?  Within \spec$_k(V)$? Are these groups like $F$ in that their subgroup spectra contain  many isomorphism classes, or does one find a generic isomorphism class of subgroup in \spec$_k(T)$ and \spec$_k(V)$?

\bibliography{refs} \bibliographystyle{plain}

\def\cprime{$'$}
\begin{thebibliography}{10}

\bibitem{MR1794135}
G.~N. Arzhantseva.
\newblock On groups in which subgroups with a fixed number of generators are
  free.
\newblock {\em Fundam. Prikl. Mat.}, 3(3):675--683, 1997.

\bibitem{MR1445193}
G.~N. Arzhantseva and A.~Yu. Ol{\cprime}shanski{\u\i}.
\newblock Generality of the class of groups in which subgroups with a lesser
  number of generators are free.
\newblock {\em Mat. Zametki}, 59(4):489--496, 638, 1996.

\bibitem{MR1929714}
Alexandre~V. Borovik, Alexei~G. Myasnikov, and Vladimir Shpilrain.
\newblock Measuring sets in infinite groups.
\newblock In {\em Computational and statistical group theory (Las Vegas,
  NV/Hoboken, NJ, 2001)}, volume 298 of {\em Contemp. Math.}, pages 21--42.
  Amer. Math. Soc., Providence, RI, 2002.

\bibitem{brin-squier}
Matthew~G. Brin and Craig~C. Squier.
\newblock Presentations, conjugacy, roots, and centralizers in groups of
  piecewise linear homeomorphisms of the real line.
\newblock {\em Comm. Algebra}, 29(10):4557--4596, 2001.

\bibitem{MR1670622}
Jos{\'e} Burillo.
\newblock Quasi-isometrically embedded subgroups of {T}hompson's group {$F$}.
\newblock {\em J. Algebra}, 212(1):65--78, 1999.

\bibitem{MR2102168}
Jos{\'e} Burillo.
\newblock Growth of positive words in {T}hompson's group {$F$}.
\newblock {\em Comm. Algebra}, 32(8):3087--3094, 2004.

\bibitem{MR1426438}
J.~W. Cannon, W.~J. Floyd, and W.~R. Parry.
\newblock Introductory notes on {R}ichard {T}hompson's groups.
\newblock {\em Enseign. Math. (2)}, 42(3-4):215--256, 1996.

\bibitem{ProbList}
Sean Cleary, John Stallings, and Jennifer Taback.
\newblock {T}hompson's group at 40 years.
\newblock {A}merican {I}nstitute of {M}athematics workshop, open problems list.
  \texttt{http://www.aimath.org/pastworkshops/thompsonsgroup.html}.

\bibitem{MR2016187}
Sean Cleary and Jennifer Taback.
\newblock Geometric quasi-isometric embeddings into {T}hompson's group {$F$}.
\newblock {\em New York J. Math.}, 9:141--148 (electronic), 2003.

\bibitem{MR909335}
Philippe Flajolet.
\newblock Analytic models and ambiguity of context-free languages.
\newblock {\em Theoret. Comput. Sci.}, 49(2-3):283--309, 1987.
\newblock Twelfth international colloquium on automata, languages and
  programming (Nafplion, 1985).

\bibitem{FlajBook}
Philippe Flajolet and Robert Sedgewick.
\newblock Analytic combinatorics.
\newblock In preparation.
  \texttt{http://algo.inria.fr/flajolet/Publications/books.html}.

\bibitem{MR2104775}
V.~S. Guba.
\newblock On the properties of the {C}ayley graph of {R}ichard {T}hompson's
  group {$F$}.
\newblock {\em Internat. J. Algebra Comput.}, 14(5-6):677--702, 2004.
\newblock International Conference on Semigroups and Groups in honor of the
  65th birthday of Prof. John Rhodes.

\bibitem{MR1725439}
V.~S. Guba and M.~V. Sapir.
\newblock On subgroups of the {R}. {T}hompson group {$F$} and other diagram
  groups.
\newblock {\em Mat. Sb.}, 190(8):3--60, 1999.

\bibitem{Jit}
Toshiaki Jitsukawa.
\newblock Stallings foldings and subgroups of free groups. {P}h{D} {T}hesis,
  {C}{U}{N}{Y} {G}raduate {C}enter, 2005.

\bibitem{MR1981427}
Ilya Kapovich, Alexei~G. Miasnikov, Paul Schupp, and Vladimir Shpilrain.
\newblock Generic-case complexity, decision problems in group theory, and
  random walks.
\newblock {\em J. Algebra}, 264(2):665--694, 2003.

\bibitem{MR2423197}
Alexei Myasnikov, Vladimir Shpilrain, and Alexander Ushakov.
\newblock Random subgroups of braid groups: an approach to cryptanalysis of a
  braid group based cryptographic protocol.
\newblock In {\em Public key cryptography---{PKC} 2006}, volume 3958 of {\em
  Lecture Notes in Comput. Sci.}, pages 302--314. Springer, Berlin, 2006.

\bibitem{MiasArx}
Alexei Myasnikov and Alexander Ushakov.
\newblock Random subgroups and analysis of the length-based and quotient
  attacks.
\newblock {\em J. Math. Crypt}, 2(1):26--61, 2008.

\bibitem{SaZi94}
Bruno Salvy and Paul Zimmermann.
\newblock Gfun: a {M}aple package for the manipulation of generating and
  holonomic functions in one variable.
\newblock {\em ACM Transactions on Mathematical Software}, 20(2):163--177,
  1994.

\bibitem{MR1676282}
Richard~P. Stanley.
\newblock {\em Enumerative combinatorics. {V}ol. 2}, volume~62 of {\em
  Cambridge Studies in Advanced Mathematics}.
\newblock Cambridge University Press, Cambridge, 1999.
\newblock With a foreword by Gian-Carlo Rota and appendix 1 by Sergey Fomin.

\bibitem{WZ1985}
J.~Wimp and D.~Zeilberger.
\newblock {Resurrecting the asymptotics of linear recurrences}.
\newblock {\em Journal of mathematical analysis and applications},
  111(1):162--176, 1985.

\bibitem{Woodruff}
Ben Woodruff.
\newblock Statistical properties of {T}hompson's group and random pseudo
  manifolds. {P}h{D} {T}hesis, {B}{Y}{U}, 2005.

\end{thebibliography}

\end{document}